\documentclass[preprint, 11pt]{elsarticle}

\usepackage[T1]{fontenc}
\usepackage[utf8]{inputenc}
\usepackage[hmargin=0.9in]{geometry}
\usepackage[babel,german=quotes]{csquotes}
\usepackage{subfig}
\usepackage{graphicx}
\usepackage[colorlinks=true,citecolor=blue,urlcolor=blue]{hyperref}
\usepackage{caption}
\usepackage{amsmath}          
\usepackage{amsthm}           
\usepackage{amssymb}          
\usepackage{bm, amsfonts}     
\usepackage{xcolor}
\usepackage{fixmath}
\usepackage{tikz}
\usepackage{bm}
\usepackage{makecell}

\newtheorem{thm}{Theorem}
\newdefinition{cor}{Corollary}
\newdefinition{example}{Example}
\newproof{pf}{Proof}
\newcommand{\proofend}{\hfill $\Box$}

\newcommand {\dx} {\,{\rm d}{\mathbf x}}

   \newcommand{\bff}{\mathbf{f}}

  \newcommand{\R}{\mathbb{R}}
  \newdefinition{rmk}{Remark}
  
  \newcommand{\pd}[2]{\frac{\partial #1}{\partial #2}}

  \newcommand{\td}[2]{\frac{\mathrm d #1}{\mathrm d #2}}

\newcommand{\beq}{\begin{equation}}
\newcommand{\eeq}{\end{equation}}

\newcommand{\bfg}{{\bf g}}

\newcommand{\bfx}{{\bf x}}
\newcommand{\sLS}{s_{h}^{\text{LS},e}}
\newcommand{\sEV}{s_{h}^{\text{EV},e}}
\newcommand{\sSUPG}{s_{h}^{\text{SUPG},e}}
\newcommand{\sVMS}{s_{h}^{\text{VMS},e}}

\newcommand{\nuEV}{\nu^{\text{EV},e}}

\newcommand{\bfq}{{\bf q}}

\newcommand{\nuSUPG}{\nu^{\text{SUPG},e}}
\newcommand{\nuVMS}{\nu^{\text{VMS},e}}

\makeatletter
\def\ps@pprintTitle{%
  \let\@oddhead\@empty
  \let\@evenhead\@empty
  \def\@oddfoot{
    \footnotesize\itshape
    \hfill\today
  }%
  \let\@evenfoot\@oddfoot}
\makeatother

\begin{document}

\begin{frontmatter} 
  \title{Entropy conservation property and entropy
    stabilization of high-order continuous
 Galerkin approximations to scalar conservation laws}

\author[TUDo]{Dmitri Kuzmin\corref{cor1}}
\ead{kuzmin@math.uni-dortmund.de}
\cortext[cor1]{Corresponding author}

\address[TUDo]{Institute of Applied Mathematics (LS III), TU Dortmund University\\ Vogelpothsweg 87,
  D-44227 Dortmund, Germany}

\author[KAUST]{Manuel Quezada de Luna}
\ead{manuel.quezada@kaust.edu.sa}

\address[KAUST]{King Abdullah University of Science and Technology (KAUST)\\ Thuwal 23955-6900, Saudi Arabia}

\journal{Computers and Fluids}

\begin{abstract}
  This paper addresses the design of linear and nonlinear stabilization
  procedures for high-order continuous Galerkin (CG) finite element
  discretizations of scalar conservation laws. We prove that the
  standard CG method is entropy conservative for the square
  entropy. In general, the rate of entropy production/dissipation
  depends on the residual of the governing equation and on the accuracy
  of the finite element approximation to the entropy variable. The
  inclusion of linear high-order stabilization generates an additional
  source/sink in the entropy budget equation. To balance the amount
  of entropy production in each cell, we construct entropy-dissipative
  element contributions using a coercive bilinear form and
  a parameter-free entropy viscosity coefficient. The 
  entropy stabilization term is high-order consistent,
  and optimal convergence behavior is achieved in practice.
  To enforce preservation of local bounds in addition to
  entropy stability, we use the Bernstein basis
  representation of the finite element solution and a new subcell
  flux limiting procedure. The underlying inequality constraints
   ensure the validity of localized entropy conditions and
  local maximum principles. The benefits of the proposed modifications
  are illustrated by numerical results for linear and nonlinear
  test problems.
\end{abstract}
\end{frontmatter}

\section{Introduction}

The design of property-preserving continuous Galerkin (CG) methods for
hyperbolic conservation laws is particularly difficult in the context of
high-order finite element approximations. Even linear advection
problems with smooth exact solutions require the use of high-order
stabilization to achieve optimal convergence rates with CG approximations
on general meshes \cite{CG-BFCT}. In the nonlinear case, a well-designed
numerical scheme should be entropy stable \cite{abgrall,chen,tadmor}.
A failure to satisfy this requirement may cause convergence to a wrong
weak solution. Tadmor's entropy stability theory \cite{tadmor87,tadmor}
provides a general framework for designing numerical fluxes that satisfy
cell entropy inequalities. Positivity preservation and local
maximum principles can be enforced using flux or slope limiting
techniques \cite{chen,kivva,entropyDG}. Many high-resolution finite volume
or discontinuous Galerkin (DG) methods are designed in this way.
Unfortunately, direct manipulation of numerical fluxes and/or solution
gradients is not an option for high-order continuous finite element
approximations. However, the desired properties can be achieved
using the framework of algebraic flux correction (AFC) for
linear transport equations \cite{afc_analysis,afc1} and its
recent extensions to nonlinear hyperbolic conservation laws
\cite{Guermond2018,Guermond2017,convex,convex2}.

An entropy stable and locally bound-preserving AFC scheme for
continuous linear ($\mathbb{P}_1$) and multilinear ($\mathbb{Q}_1$)
finite elements was designed in \cite{entropyCG} using graph
Laplacian stabilization and a monolithic limiting strategy.
Alternative approaches to enforcing entropy stability in
finite element schemes include the use of residual-based
entropy viscosity \cite{Guermond2018,Guermond2017} and Rusanov-type
penalization for the gradients of entropy variables
\cite{abgrall,ranocha,entropyDG}. In the present paper,
we extend the entropy correction tools
proposed in \cite{abgrall,entropyDG,entropyCG}
to stabilized high-order CG approximations and combine them with the
subcell flux limiting strategy developed in \cite{convex2}
for AFC schemes based on high-order Bernstein
finite elements. Moreover, we prove that the standard CG
method is entropy conservative for the square entropy. In
contrast to the square entropy stability property of DG
methods for scalar conservation laws \cite{jiang},
this result seems to be largely unknown. The use of a
general entropy and/or inclusion of linear high-order stabilization terms
produces additional sources or sinks in the entropy balance equation
associated with the semi-discrete CG scheme. 
To convert this equation into a discrete entropy inequality,
we add a nonlinear entropy dissipation term which represents
a generalized high-order version of Abgrall's \cite{abgrall}
entropy fix. In the process of subcell flux correction,
we blend the resulting entropy stable high-order scheme and
a low-order compact-stencil
approximation of Rusanov (local Lax-Friedrichs) type
in a manner which guarantees the validity of
all relevant constraints (conservation principles,
entropy inequalities, maximum principles). Numerical studies
are performed for linear and nonlinear test problems.

\section{Entropy conservation property of the CG method}
\label{sec:cg_ent_stab}

Let $u(\mathbf{x},t)$ be a scalar conserved
quantity depending on the space location $\mathbf{x}\in
\R^d,\ d\in\{1,2,3\}$ and time instant $t\ge 0$. Consider
an initial value problem of the form 
\begin{subequations}
\begin{align}
 \pd{u}{t}+\nabla\cdot\mathbf{f}(u)=0 &\qquad\mbox{in}\ \R^d\times\R_+,
\label{ibvp-pde}\\
 u(\cdot,0)=u_0 &\qquad\mbox{in}\ \R^d,\label{ibvp-ic}
\end{align}
\end{subequations}
where $\mathbf{f}=(\mathsf{f}_1,\ldots,\mathsf{f}_d)$ is a possibly
nonlinear flux function and $u_0:\R^d\to\mathcal G$ is the data of
the initial condition. A convex set $\mathcal G\subset\R$ is
called an {\it invariant set} of problem  \eqref{ibvp-pde}--\eqref{ibvp-ic}
if the exact solution $u$ stays in $\mathcal G$ for all $t> 0$ \cite{Guermond2016}.
A convex function
$\eta:\mathcal G\to\R$ is called an {\it entropy} and $v(u)=\eta'(u)$ is called an
{\it entropy}
variable if there exists an {\it entropy flux} 
$\mathbf{q}:\mathcal G\to\mathbb{R}^d$ such that $v(u)\mathbf{f}'(u)=\mathbf{q}'(u)$.
 A~weak solution
$u$ of \eqref{ibvp-pde} is called an {\it entropy solution}
if the entropy inequality
\beq
 \pd{\eta}{t}+\nabla\cdot\mathbf{q}(u)\le 0 \qquad\mbox{in}\ \R^d\times\R^+
\label{ent-ineq}
\eeq
holds for any {\it entropy pair} $(\eta,\mathbf{q})$. For any smooth weak
solution, the conservation law
\beq
 \pd{\eta}{t}+\nabla\cdot\mathbf{q}(u)=0 \qquad\mbox{in}\ \R^d\times\R^+
 \eeq
 can be derived from \eqref{ibvp-pde} using multiplication by the entropy
 variable $v$, the chain rule, and the definition of an entropy pair. Hence,
 entropy is conserved in smooth regions and dissipated at shocks.

Adopting the terminology of
Guermond et al. \cite{Guermond2018,Guermond2016}, we will call 
a numerical scheme {\it invariant
domain preserving} (IDP) if the solution of the (semi-)discrete problem
is guaranteed to stay in an invariant set $\mathcal G$. Additionally, a property-preserving
discretization of \eqref{ibvp-pde} should be {\it entropy stable}, i.e., it should satisfy
a discrete version of the entropy inequality \eqref{ent-ineq}. The lack of
entropy stability is a typical reason for convergence of numerical schemes
to nonphysical weak solutions.

To discretize \eqref{ibvp-pde} in a bounded domain $\Omega\subset\R^d$,
we use the continuous Galerkin (CG) method and a conforming mesh
$\mathcal T_h=\{K^1,\ldots,K^{E_h}\}$.
For simplicity, we assume
that the imposed boundary conditions are periodic.
Let $V_{h}^e\in\{\mathbb{P}_p(K^e),
\mathbb{Q}_p(K^e)\},\ p\in\mathbb{N}$ denote the polynomial space of
the finite element approximation on $K^e,\ e=1,\ldots,E_h$ and
$V_h=\{v_h\in C(\bar\Omega_h)
\,:\, v_h|_{K^e}\in V_h^e,\  e=1,\ldots,E_h\}$ the space of
continuous piecewise-polynomial functions defined
on $\bar\Omega_h=\bigcup_{i=1}^{E_h}K^e$. Each function $v_h\in V_h$ can
be written as $v_h=\sum_{j=1}^{N_h}v_j\varphi_j$, where $\varphi_1,\ldots,
\varphi_{N_h}$ are Lagrange or Bernstein basis functions associated with
nodal points $\mathbf{x}_1,\ldots,\mathbf{x}_{N_h}\in\bar\Omega_h$.
We define the {\it full stencil} of node $i$
as the integer set $\mathcal N_i=\bigcup_{e\in\mathcal E_i}\mathcal N^e$,
where $\mathcal E_i$
denotes the set of (numbers of) elements containing the point
$\mathbf{x}_i$ and $\mathcal N^e$ is the set of (numbers of)
nodes belonging to $K^e$. In addition to full stencils,
we will use compact nearest-neighbor stencils 
in the description of the proposed methods below.

Approximating the exact entropy
solution $u$ of \eqref{ibvp-pde} by $u_h\in V_h$,
we consider the CG discretization
\beq\label{CGhigh}
\sum_{e=1}^{E_h}\int_{K^e} w_h\left(\frac{\partial u_h}{\partial t}+\nabla\cdot \bff(u_h)\right)\dx=0
  \qquad \forall w_h\in V_h.
\eeq
Non-periodic flux boundary conditions can be taken into account
by adding an integral
over the inflow boundary of the computational
domain $\Omega_h$, see \cite{convex,convex2} for details.

\begin{thm}[Entropy behavior of the continuous Galerkin method]\label{thm1}
  Let $\{\eta(u),\mathbf{q}(u)\}$ be an entropy pair and $v(u)=\eta'(u)$
  the corresponding entropy variable.
  Suppose that \eqref{CGhigh}
  holds for $u_h\in V_h$. Then
  \beq\label{cons-eta}
 \sum_{e=1}^{E_h}\int_{K^e}\left(\pd{\eta(u_h)}{t}+
\nabla\cdot\mathbf{q}(u_h)\right)\dx=\sum_{e=1}^{E_h}
\int_{K^e}(v(u_h)-v_h)\left(\pd{u_h}{t}+\nabla\cdot\mathbf{f}(u_h)\right)\dx
\quad\forall v_h\in V_h.
\eeq
In particular, $\eta(u_h)$ satisfies \eqref{cons-eta}
 for $v_h\equiv 0$ and for an arbitrary
approximation $v_h\in V_h$ to $v(u_h)$.

\end{thm} 

\begin{pf}
  Using the chain rule to differentiate $\eta(u)$, substituting
  $v(u_h)$ for $\eta'(u_h)$ and recalling that $v(u_h)\mathbf{f}'(u_h)
  =\mathbf{q}'(u_h)$ by definition of the entropy flux,
  we transform the left-hand side of \eqref{cons-eta} as follows:
  \begin{align*}
    \sum_{e=1}^{E_h}\int_{K^e}\left(\pd{\eta(u_h)}{t}+
    \nabla\cdot\mathbf{q}(u_h)\right)\dx&=
    \sum_{e=1}^{E_h}\int_{K^e}\left(v(u_h)\pd{u_h}{t}+
    \mathbf{q}'(u_h)\cdot\nabla u_h\right)\dx\\
    &= \sum_{e=1}^{E_h}\int_{K^e}v(u_h)\left(\pd{u_h}{t}+
    \mathbf{f}'(u_h)\cdot\nabla u_h\right)\dx\\
    &= \sum_{e=1}^{E_h}\int_{K^e}v(u_h)\left(\pd{u_h}{t}+
     \nabla\cdot\mathbf{f}(u_h)\right)\dx.
  \end{align*}
  The validity of representation \eqref{cons-eta} follows from the fact that
  \eqref{CGhigh} holds for any $v_h\in V_h$. \hfill
\proofend
\end{pf}

Theorem \ref{thm1} reveals that the CG method is globally {\it entropy
  conservative} in the following sense.

\begin{cor}[Entropy conservation property of the continuous
  Galerkin method]
  For the square entropy $\eta(u)=\frac{u^2}{2}$, the CG approximation
  $u_h$ satisfies
    \beq
 \sum_{e=1}^{E_h}\int_{K^e}\left(\pd{\eta(u_h)}{t}+
\nabla\cdot\mathbf{q}(u_h)\right)\dx=0.
\eeq
\end{cor}

\begin{pf}
  The entropy variable associated with $u_h\in V_h$ is
  $v(u_h)=u_h$. Hence, the right-hand side of the entropy balance equation
  \eqref{cons-eta} vanishes for the CG solution $u_h$ defined by \eqref{CGhigh}.
  \end{pf}

\begin{rmk}
  The rate of net entropy production for $\eta(u)\ne\frac{u^2}{2}$ is given
  by the weighted residual $$\sum_{e=1}^{E_h}
  \int_{K^e}v(u_h)\left(\pd{u_h}{t}+\nabla\cdot\mathbf{f}(u_h)\right)\dx.$$
  A similar result was obtained by Jiang and Shu \cite{jiang} in the context
  of DG methods.
\end{rmk}

\section{Linear high-order stabilization}\label{sec:ls_stab}

It is common knowledge that the convergence behavior of the CG 
method is unsatisfactory even for linear advection problems with smooth exact
solutions. The provable order of accuracy w.r.t. the $L^2$ error is $O(h^p)$,
where $h$ is the mesh size and $p$ is the polynomial degree of the finite
element approximation (see, e.g., Quarteroni and Valli \cite{QV94},
eq.~(14.3.16)). To achieve optimal $O(h^{p+1/2})$ convergence rates on
general meshes \cite{Burman2010,John2011}, many stabilized CG methods
of the form
\beq\label{CGstab}
\sum_{e=1}^{E_h}\left[\int_{K^e} w_h\left(\frac{\partial u_h}{\partial t}+\nabla\cdot \bff(u_h)\right)\dx+\sLS(w_h,u_h)\right]=0
  \qquad \forall w_h\in V_h
\eeq
were proposed in the literature. The linear
stabilization term $\sLS(w_h,u_h)$ is
usually defined as a weighted residual of the governing equation or of
another relation which is satisfied for $h\to 0$.

The stabilization term of the streamline upwind Petrov-Galerkin (SUPG) method 
\cite{supg} is defined by
\beq\label{shSUPG}
\sSUPG (w_h,u_h)=
\nuSUPG\int_{K^e}\left(\bff^\prime(u_h)\cdot\nabla w_h\right)
\left(\dot u_h+\nabla\cdot \bff(u_h)\right)\dx,
\eeq
where $\dot u_h\in V_h$ is an approximate time derivative and
$\nuSUPG$ is a stabilization parameter depending on the local mesh size
$h^e=|K^e|^{1/d}$. In the numerical experiments of Section \ref{sec:num},
we use 
\beq\label{nuSUPG}
\nuSUPG=\frac{\omega h^e}{2p\|\bff^\prime(u_h)\|_{L^\infty(K^e)}},
\eeq
where  $\omega=1$ by default. Smaller values of $\omega$ can be
used to adjust the amount of linear stabilization.

Existing theory  \cite{Burman2010,John2011} guarantees
$O(h^{p+1/2})$ convergence behavior of the consistent SUPG method for
the linear advection equation provided that $\dot u_h$ is a sufficiently
good approximation to $\pd{u_h}{t}$. The coefficients of $\dot u_h$
corresponding
to \eqref{CGhigh} are given by the solution of the linear system
\beq\label{CGudot}
\sum_{e\in\mathcal E_i}\sum_{j\in\mathcal N^e}m_{ij}^e\dot u_j
=-\sum_{e\in\mathcal E_i}
\int_{K^e}\varphi_i\nabla\cdot \bff(u_h)\dx,
  \qquad i=1,\ldots,N_h,
\eeq
where we use the stencil notation introduced in Section \ref{sec:cg_ent_stab}.
The coefficients $m_{ij}^e$ are defined by
\beq\label{mijdef}
m_{ij}^e=\int_{K^e}\varphi_i\varphi_j\dx.
\eeq
In contrast to other approximations of the time derivative, definition
\eqref{CGudot}
supports the use of general time integrators and avoids the dependence
of the parameter $\nuSUPG$ on the time step.

A closely related {\it variational multiscale} (VMS) method
\cite{John2006,CG-BFCT} stabilizes \eqref{CGhigh} using the
bilinear form
\beq\label{shVMS}
\sVMS(w_h,u_h)=\nuVMS
\int_{K^e}\nabla w_h\cdot(\nabla u_h-\bfg_h)\dx,
\eeq
where $\mathbf{g}_h=(\mathsf{g}_{h1},\ldots,\mathsf{g}_{hd})\in
(V_h)^d$ is a continuous approximation to the gradient
 $\nabla u_h$ and 
\beq\label{nuVMS}
\nuVMS = \frac{\omega h^e\|\bff^\prime(u_h)\|_{L^\infty(K^e)}}{2p}.
\eeq
 
Lohmann et al. \cite{CG-BFCT} defined the continuous gradient
$\mathbf{g}_h$ using the $L^2$ projection 
\begin{align}\label{cons_grad_rec}
 \sum_{e=1}^{E_h} \int_{K^e} w_h(\bfg_h-\nabla u_h)\dx = 0 \qquad \forall w_h\in V_h
\end{align}
which requires solution of $d$ linear systems with the consistent
mass matrix
\beq
M_C=(m_{ij})_{i,j=1}^{N_h},\qquad m_{ij}=\sum_{e\in\mathcal E_i
  \cap\mathcal E_j}m_{ij}^e.
\eeq
As shown in \cite{CG-BFCT}, the one-dimensional version of
\eqref{shVMS} using this definition of $\bfg_h$ is equivalent
to the SUPG stabilization
\eqref{shSUPG} for linear advection with constant velocity.
  
As an inexpensive alternative to gradient recovery via
consistent-mass $L^2$ projections, we define 
\beq\label{averaged_grad_rec}
\bfg_h=\sum_{j=1}^{N_h} \bfg_j\psi_j
\eeq
using
Lagrange basis functions $\psi_1,\ldots,\psi_{N_h}$to interpolate the
averaged nodal values
\beq\label{mlump}
\bfg_i =\frac{1}{m_i}\sum_{e\in\mathcal{E}_i}m_i^e \nabla u_h|_{K^e}(\bfx_i),\qquad
m_i=\sum_{e\in\mathcal{E}_i}m_i^e,\qquad m_i^e=\int_{K^e}\varphi_i^e\dx.
\eeq
Note that both definitions of the reconstructed gradient produce
$\bfg_h=\nabla u_h$ in the case $\nabla u_h\in (V_h)^d$.

\begin{rmk}
If $u_h\in V_h$
is defined using the Lagrange basis as well, then $\varphi_j
=\psi_j$ for $j=1,\ldots,N_h$. The limiting techniques
presented in Section \ref{sec:ev_stab} require the use of
Bernstein basis functions $\varphi_j\ne\psi_j$.
\end{rmk}

The ability of the above stabilization techniques to deliver the expected
convergence rates for high-order finite element approximations to linear
advection problems is verified in Section \ref{sec:num}.

\section{Nonlinear high-order stabilization}\label{sec:ev_stab}

By virtue of Theorem \ref{thm1}, the entropy of the stabilized finite
element approximation $u_h$ satisfies
  \beq\label{cons-eta-stab}
 \sum_{e=1}^{E_h}\int_{K^e}\left(\pd{\eta(u_h)}{t}+
\nabla\cdot\mathbf{q}(u_h)\right)\dx=\sum_{e=1}^{E_h}p^e_h(v_h,u_h),
\eeq
where $v_h=\sum_{j=1}^{N_h}v_j\varphi_j$ is defined in terms of
$v_j=v(u_j),\ j=1,\ldots,N_h$ or as $L^2$ projection of $v(u_h)$ into $V_h$. The
bilinear form $p^e_h(v_h,u_h)$ of the local entropy production term is defined by
\beq\label{prodEV}
p^e_h(v_h,u_h)=\int_{K^e}(v(u_h)-v_h)\left(\pd{u_h}{t}+\nabla\cdot\mathbf{f}(u_h)\right)\dx-\sLS(v_h,u_h).
\eeq
If we have $p^e_h(v_h,u_h)\le 0$ for all $e=1,\ldots,E_h$, then the discrete
entropy inequality
\beq\label{eta-ineq}
\sum_{e=1}^{E_h}\int_{K^e}\left(\pd{\eta(u_h)}{t}+
\nabla\cdot\mathbf{q}(u_h)\right)\dx\le 0
\eeq
holds for the solution $u_h$ of the semi-discrete problem
\eqref{CGstab}. Ironically, the contribution of the
linear stabilization term $\sLS (v_h,u_h)$ may render $p^e_h(v_h,u_h)$
positive even in the case of the square entropy
$\eta(u)=\frac{u^2}{2}$, in which $v(u_h)=u_h$ and
\eqref{eta-ineq} holds as equality for the solution of \eqref{CGhigh}.
In other words, the use of (nonsymmetric) linear stabilization can
cause or aggravate the lack of entropy stability.

To limit the amount of entropy production, we introduce 
an entropy viscosity (EV) term $\sEV(w_h,v_h)$ such
that $p^e_h(v_h,u_h)-\sEV(v_h,v_h)\le 0$ and, therefore, 
\eqref{eta-ineq} holds for the solution $u_h$ of 
\beq\label{CGstabEV}
\sum_{e=1}^{E_h}\left[\int_{K^e} w_h\left(\frac{\partial u_h}{\partial t}+\nabla\cdot \bff(u_h)\right)\dx+\sLS(w_h,u_h)+\sEV(w_h,v_h)\right]=0
  \qquad \forall w_h\in V_h.
  \eeq
 Entropy correction techniques of
 this kind trace their origins to the work of Abgrall \cite{abgrall}.
As pointed out in \cite{entropyDG}, any symmetric positive definite
(coercive) bilinear form $b^e(\cdot,\cdot)$ can be used to construct
$\sEV(w_h,v_h)$. Following
\cite{ranocha,entropyCG}, we choose $b^e(\cdot,\cdot)$
to be the $L^2(K^e)$ scalar product and define
\beq\label{stabEV}
\sEV(w_h,v_h)=\nuEV
\int_{K^e}(i_{h,1}^{e}w_h-i_{h,0}^{e}w_h)
(i_{h,1}^{e}v_h-i_{h,0}^{e}v_h)\dx,
\eeq
where $i_{h,1}^{e}v_h$ is the piecewise $\mathbb{P}_1/\mathbb{Q}_1$
Lagrange
interpolant of the nodal values $\{v_h(\mathbf{x}_j),\ j\in\mathcal N^e\}$ 
and $i_{h,0}^{e}v_h$ is a piecewise-constant approximation
defined by the subcell averages of $i_{h,1}^{e}v_h$, i.e.,
by averages over the elements of the submesh formed by the
nodes $\{\mathbf{x}_j,\ j\in\mathcal N^e\}$, cf.~\cite{convex2}.

    \begin{rmk}
  The entropy stabilization term \eqref{stabEV} and the entropy
  viscosity coefficient $\nuEV$ can also be constructed using
  $i_{h,1}^{e}v_h:=v_h$ and
  $i_{h,0}^{e}v_h:=\frac{1}{|K^e|}
  \int_{K^e}v_h\dx$. This definition corresponds to the Rusanov
  dissipation employed in \cite{ranocha,entropyDG}. The
  two versions are equivalent for $\mathbb{P}_1/\mathbb{Q}_1$
   elements.
\end{rmk}

    It remains to define the EV parameter $\nuEV\ge 0$. A lower bound
    $\nu^{\text{EV},e,\min}$
which guarantees entropy stability
of the semi-discrete scheme \eqref{CGstabEV} is provided by the
following theorem.

\begin{thm}[Minimal entropy viscosity]\label{thm2}
  Let $p_h^e(v_h,u_h)$ be the local entropy production term defined by
  \eqref{prodEV} and  $\nu^{\text{\rm EV},e}$ an entropy viscosity
  coefficient which is bounded below by
\beq\label{nu_min}
  \nu^{\text{\rm EV},e,\min} := 
  \frac{\max\{0,p_h^e(v_h,u_h)\}}{\int_{K^e}(i_{h,1}^{e}v_h-i_{h,0}^{e}v_h)^2\dx}.
  \eeq
   Then the semi-discrete scheme defined by \eqref{CGstabEV}
  and \eqref{stabEV}
  satisfies the discrete entropy inequality \eqref{eta-ineq}.
\end{thm}
 
\begin{pf}
  The assertion of the theorem is a direct consequence of
  the fact that the local entropy condition 
 \beq\label{stabEVloc}
 p_h^e(v_h,u_h)-
 \sEV(v_h,v_h)\le 0
 \eeq
 holds for  $\sEV(v_h,v_h)
 \ge \nu^{\text{\rm EV},e,\min}\int_{K^e}(i_{h,1}^{e}v_h-i_{h,0}^{e}v_h)^2\dx
 = \max\{0,p_h^e(v_h,u_h)\}.$
 \hfill\proofend
 
 \end{pf}

The use of $\nuEV= \nu^{\text{EV},e,\min}$ in \eqref{stabEV} introduces the
minimal amount of entropy stabilization which ensures the validity
of the discrete entropy inequality \eqref{eta-ineq}. The corresponding
semi-discrete scheme \eqref{CGstabEV} is {\it barely} entropy stable
(cf. \cite{entropyDG,entropyCG}) and may fail to converge to correct
weak solutions. Adopting Tadmor's design philosophy \cite{tadmor87,tadmor},
we adjust the levels of entropy dissipation by using a stabilization
parameter $\nuEV\ge \nu^{\text{EV},e,\min}$
which depends on the local smoothness of
the approximate solution.

Let $\pi^e_h:L^2(K^e)\to\mathbb{P}_{p-1}(K^e)$ denote the local $L^2$
projection operator into the polynomial space of degree $p-1\ge 0$.
For any $u\in L^2(K^e)$, the polynomial $\pi^e_hu\in\mathbb{P}_{p-1}(K^e)$
is defined by
\begin{align}
  \int_{K^e}w_h^e(u-\pi^e_hu) \dx=0 \qquad
  \forall w_h^e\in \mathbb{P}_{p-1}(K^e).
\end{align}
To gain better control of local
entropy production without losing high-order accuracy
in smooth regions, we define the nonlinear stabilization term
\eqref{stabEV} using the entropy viscosity coefficient
\begin{align}\label{nu}
  \nuEV=
  \nu^{\text{EV},e,\min} + 
 \frac{\left|\int_{K^e}\nabla v_h
         \cdot(\bff(\pi^e_hu_h)-\bff(u_h))\dx\right|}{
   \int_{K^e}(i_{h,1}^{e}v_h-i_{h,0}^{e}v_h)^2\dx}.
\end{align}

\begin{rmk}
In the unlikely case that $\nuEV$ defined by \eqref{nu} becomes
very large, explicit time discretizations of \eqref{CGstabEV}
may require the use of impractically small time steps. If implicit
treatment of \eqref{stabEV}, as proposed in \cite{entropyDG} in the
context of DG-$\mathbb{P}_1$ approximations, is not an option, then an
upper bound $\nu^{\text{EV},e,\max}$ may need to be imposed on the value
of $\nuEV$. For $\nu^{\text{EV},e,\max}<\nu^{\text{EV},e,\min}$,
condition \eqref{stabEVloc}
cannot be satisfied using \eqref{stabEV}
with $\nuEV\le \nu^{\text{EV},e,\max}$. However, entropy
stability is still guaranteed if~\eqref{CGstabEV} is
constrained using the convex limiting techniques that we present
in the next section.
\end{rmk}

\section{Monolithic convex limiting}\label{sec:ES-lim}

The stabilized high-order finite element scheme \eqref{CGstabEV}
may require additional modifications to ensure the invariant domain
preservation (IDP) property and validity of local maximum principles
for problems with discontinuities and propagating fronts.
Bound-preserving convex limiting techniques for high-order Bernstein
finite element approximations to scalar hyperbolic conservation laws
were developed in \cite{convex2} without taking entropy conditions
into account. Entropy stability preserving (ESP) limiters were
introduced in \cite{entropyDG,entropyCG} in the context of 
$\mathbb{P}_1/\mathbb{Q}_1$ approximations of CG and DG type.
In this section, we generalize the convex limiting tools developed in
\cite{entropyDG,convex2,entropyCG} and apply them to \eqref{CGstabEV}.

Let $\varphi_1,\ldots,\varphi_{N_h}$ be Bernstein basis functions
spanning the finite element space $V_h$. The definition of these
basis functions for simplex and tensor-product meshes can be
found, e.g., in \cite{CG-BFCT} and in the Appendix of \cite{convex2}.
The corresponding degrees of freedom $u_1,\ldots,u_{N_h}$ are
associated with the nodal points $\bfx_1,\ldots,\bfx_{N_h}$ and
called {\it Bernstein coefficients}. The approximate solution
$u_h\in V_h$ satisfies \cite{CG-BFCT}
\beq
\min_{j\in\mathcal N^e}u_j\le u_h(\mathbf{x})=\sum_{j\in\mathcal N^e}
u_j\varphi_j(\mathbf{x})\le \max_{j\in\mathcal N^e}u_j\qquad
\forall\mathbf{x}\in K^e.
\eeq
Hence, the IDP property 
is guaranteed if all Bernstein coefficients  of $u_h$
are in the admissible range.

Substituting test functions $w_h\in\{\varphi_1,\ldots,\varphi_{N_h}\}$
into \eqref{CGstabEV}, we obtain a system of semi-discrete equations
for the (generally time-dependent) Bernstein coefficients. This
system is given by
\begin{align}\label{target}
 \sum_{e\in\mathcal{E}_i}\sum_{j\in\mathcal{N}^e}m_{ij}^e\frac{\text{d} u_j}{\text{d}t}=
-\sum_{e\in\mathcal{E}_i}\left[\int_{K^e}\varphi_i\nabla\cdot \bff(u_h)\dx
 +\sLS(\varphi_i,u_h)+\sEV(\varphi_i,v_h)\right],\quad i=1,\ldots,N_h.
\end{align}
The coefficients $m_{ij}^e$ of the consistent element mass matrix
$M_C^e=(m_{ij}^e)_{i,j=1}^{N_h}$ are defined by \eqref{mijdef}.

In the process of monolithic
convex limiting \cite{convex,convex2}, the element contributions
to the residual of the high-order
{\it target scheme} \eqref{target} are modified
 to guarantee the validity of property-preserving inequality
constraints. For that purpose, we introduce the lumped element
mass matrix $M_L^e=(\delta_{ij}m_i^e)_{i,j=1}^{N_h}$, the element
matrix $\mathbf{C}^e=(\mathbf{c}_{ij}^e)_{i,j=1}^{N_h}$ of the
discrete gradient operator, its  `lumped'
counterpart $\tilde{\mathbf{C}}^e=M_L^e(M_C^e)^{-1}\mathbf{C}^e$,
and a discrete diffusion operator $\tilde D^e=(\tilde d_{ij}^e)_{i,j=1}^{N_h}$.
The diagonal entries $m_i^e=\sum_{j\in\mathcal N^e}m_{ij}^e$
of $M_L^e$ are the weights that we
used in \eqref{mlump}. They are positive since the Bernstein
basis functions are nonnegative by definition. The vector-valued
entries of $\mathbf{C}^e$ are given by
\beq
\mathbf{c}_{ij}^e=\int_{K^e}\varphi_i\nabla\varphi_j\dx.
\eeq
If node $i$ or node $j$ is an interior point of $\Omega_h$, integration by
parts using Green's formula yields
\beq\label{skew}
\sum_{e=1}^{E_h}\mathbf{c}_{ji}^e=-\sum_{e=1}^{E_h}\mathbf{c}_{ij}^e.
\eeq
An analytical formula for the entries of $\tilde{\mathbf{C}}^e$ is
derived in the Appendix of \cite{convex2}, where we show that this
element matrix has the same compact sparsity pattern as the
collocated piecewise
$\mathbb{P}_1/\mathbb{Q}_1$ approximation
on a subdivision of $K^e$ into subcells.
That is, we have $\tilde{\mathbf{c}}_{ij}=0$
if nodes $i$ and~$j$ are not nearest neighbors belonging to the
same subcell. On meshes consisting of parallelograms or
parallelepipeds, the entries $\tilde{\mathbf{c}}_{ij}$ associated
with diagonal subcell neighbors vanish as well \cite{hennes}.
Therefore, the sparsity pattern of the element contribution
$\tilde{\mathbf{C}}^e$ to the lumped discrete gradient is
defined by tensor products of one-dimensional three-point
stencils for $d$ coordinate directions. The
element matrix $\tilde D^e$ of the discrete diffusion
(alias {\it graph Laplacian} \cite{Guermond2018,Guermond2016}) operator
has the same sparsity pattern as  $\tilde{\mathbf{C}}^e$.
Let $\tilde{\mathcal N}_i^e
  =\{j\in\mathcal N^e\,:\,|\tilde{\mathbf{c}}_{ij}^e|+|\tilde{\mathbf{c}}_{ji}^e|
 > 0\}$ denote the compact element stencil of node $i\in\mathcal N^e$.
By default, the artificial diffusion coefficients $\tilde d_{ij}^e$
are defined by the local Lax-Friedrichs
(LLF) formula \cite{Guermond2016,convex2,entropyCG}
\beq\label{dijmax}
\tilde d_{ij}^{e}=\begin{cases}
\max\{|\tilde{\mathbf{c}}_{ij}^e|
\lambda_{ij}^{\max},|\tilde{\mathbf{c}}_{ji}^e|\lambda_{ji}^{\max}\} & \mbox{if}\
i\in\mathcal N^e,\ j\in \tilde{\mathcal N_i^e}\backslash\{i\},\\
-\sum_{k\in {\tilde{\mathcal N}_i^e}\backslash\{i\}}\tilde d_{ik}^e & \mbox{if}\ j=i\in\mathcal N^e,\\
0 & \mbox{otherwise},
\end{cases}
\eeq
where $\lambda_{ij}^{\max}$ is a guaranteed upper bound for the
maximal wave speed \cite{Guermond2016,entropyCG}
\beq
\lambda_{ij}^{\max}\ge\max_{\omega\in[0,1]}|\mathbf{n}_{ij}^e\cdot
\mathbf{f}'(\omega u_i+(1-\omega) u_j)|,\qquad \mathbf{n}_{ij}^e=
\frac{\tilde{\mathbf{c}}_{ij}^e}{|\tilde{\mathbf{c}}_{ij}^e|}.
\eeq
The monolithic convex limiting approach developed in
\cite{convex,convex2,entropyCG} approximates \eqref{target}
by
\beq\label{fcorr}
\sum_{e\in\mathcal{E}_i}m_{i}^e\td{u_i}{t}=
\sum_{e\in\mathcal E_i}\sum_{j\in\tilde{\mathcal N}_i^e\backslash\{i\}}
[\tilde d_{ij}^e(u_j-u_i)+\bar f_{ij}^{e}-\tilde{\mathbf{c}}_{ij}^e\cdot
  (\mathbf{f}(u_j)-\mathbf{f}(u_i))].
  \eeq
As shown in \cite{Guermond2016}, the low-order LLF scheme corresponding
to $\bar f_{ij}^{e}=0$ is locally bound-preserving and entropy stable for
any convex entropy. If discretization in time is performed using a
strong stability preserving (SSP) Runge-Kutta method \cite{ssprev},
the forward Euler update corresponding to a single stage satisfies
a discrete entropy inequality for any convex entropy $\eta(u)$
\cite{Guermond2016}. Moreover, a local discrete maximum principle
holds for time steps satisfying a CFL-like time step restriction
\cite{Guermond2016,convex2}.

The stabilized high-order approximation \eqref{target} can also 
be written in the compact-stencil form \eqref{fcorr} using an array
of {\it antidiffusive} subcell fluxes $\tilde f_{ij}^{e}$ which we derive here for
the reader's convenience. Let $\tilde \varphi_1,\ldots,\tilde\varphi_{N_h}$
denote the basis functions of the subcell $\mathbb{P}_1/\mathbb{Q}_1$
approximation. Define \cite{hennes}
\beq
\tilde m_{ij}^e=\begin{cases}
\int_{K^e}\tilde \varphi_i\tilde\varphi_j\dx & \mbox{if}\
i\in\mathcal N^e,\ j\in\tilde{\mathcal N}_i^e\backslash\{i\},\\
-\sum_{k\in \tilde{\mathcal N}_i\backslash\{i\}}\tilde m_{ik}^e
&  \mbox{if}\ j=i\in \mathcal N^e,\\
0 &  \mbox{otherwise}.
\end{cases}
\eeq

The compact-stencil scheme \eqref{fcorr} with
subcell fluxes $\bar f_{ij}^{e}=\tilde f_{ij}^{e}$
is equivalent to \eqref{target} for
\beq\label{subcell}
\tilde f_{ij}^{e}=\tilde m_{ij}^e(\dot w_j^e-\dot w_i^e)+
\tilde d_{ij}^e(u_i-u_j),\qquad \forall i\in\mathcal N^e,\ j\in
\tilde{\mathcal N}^e_i\backslash\{i\},
\eeq
where $\dot w_i^e,\ i\in \mathcal N^e$ are subcell flux potentials
satisfying the small sparse linear system (cf. \cite{convex2})
\begin{align}
\sum_{j\in\tilde{\mathcal{N}}_i^e}\tilde m_{ij}^e\dot w_j^e&=
\sum_{j\in\mathcal N^e\backslash\{i\}}
  m_{ij}^e(\dot u_i^S-\dot u_j^S)
  +\sum_{j\in\mathcal N^e}
  (\tilde{\mathbf{c}}_{ij}^e-\mathbf{c}_{ij}^e)\cdot\mathbf{f}(u_j)
  -\sum_{j\in \mathcal N^e}\mathbf{c}_{ji}^e\cdot\mathbf{f}(u_j)
\nonumber\\
&+ 
\int_{K^e}\nabla\varphi_i\cdot\mathbf{f}(u_h)\dx
-\sLS(\varphi_i,u_h)-\sEV(\varphi_i,v_h),\qquad
 i\in\mathcal N^e. \label{aux}
\end{align}
  The number of unknowns equals the number of nodes per element.
The solution of  \eqref{aux} is  determined up to a constant,
whose value has no influence on the value of  the flux
$\tilde f_{ij}^{e}$ defined by \eqref{subcell}.
The Bernstein coefficients $\dot u_j^S$ of the stabilized
approximate time derivatives are obtained by solving
\begin{align}\label{targetudot}
  \sum_{e\in\mathcal{E}_i}\sum_{j\in\mathcal{N}^e}m_{ij}^e\dot u_j^S=
-\sum_{e\in\mathcal{E}_i}\left[\int_{K^e}\varphi_i\nabla\cdot \bff(u_h)\dx
 +\sLS(\varphi_i,u_h)+\sEV(\varphi_i,v_h)\right],\quad i=1,\ldots,N_h.
\end{align}

The subcell flux limiter proposed in \cite{convex2} constrains
$\tilde f_{ij}^{e}$ in a manner which guarantees that the
result $\bar u_i$ of each SSP Runge-Kutta stage is bounded
by the input data $u_j,\ j\in\mathcal N_i$ as follows:
\beq\label{afc-bounds}
\min_{j\in \mathcal N_i}u_j=:u_i^{\min}
\le \bar  u_i\le u_i^{\max}:=\max_{j\in \mathcal N_i}u_j.
\eeq
A proof of the IDP property  is based on the representation of
$\bar  u_i$  in terms of the {\it bar states} 
\beq\label{bar-states}
\bar u_{ij}^e
=\frac{u_j+u_i}{2}-\frac{\tilde{\mathbf{c}}_{ij}^e
  \cdot(\mathbf{f}(u_j)-\mathbf{f}(u_i))}{2\tilde d_{ij}^{e}}
\in[u_i^{\min},u_i^{\max}],\qquad
\bar u_{ij}^{e,*}=\bar u_{ij}^e+\frac{\alpha_{ij}^e
  \tilde f_{ij}^{e,*}}{2\tilde d_{ij}^{e}}
\eeq
such that $\bar u_{ij}^{e,*}\in[u_i^{\min},u_i^{\max}]$ for any
$\alpha_{ij}^e\in[0,1]$ if the limited flux $\tilde f_{ij}^{e,*}$ is given by
\cite{convex,convex2}
\beq\label{fij_lim}
\tilde f_{ij}^{e,*}=\begin{cases}
  \min\,\left\{\tilde f_{ij}^e,2\tilde d_{ij}^e\min\,\{u_i^{\max}-\bar u_{ij}^e,
    \bar u_{ji}^e-u_j^{\min}\}\right\} & \mbox{if}\  \tilde f_{ij}^e>0,\\[0.25cm]
    \max\left\{\tilde f_{ij}^e,2\tilde d_{ij}^e\max\{u_i^{\min}-\bar u_{ij}^e,
    \bar u_{ji}^e-u_j^{\max}\}\right\} & \mbox{otherwise}.
\end{cases}
\eeq
For further explanations and detailed proofs,
we refer the interested reader
to \cite{convex,convex2,entropyCG}. After the application
of the IDP limiter, the magnitude of the bound-preserving
flux $\tilde f_{ij}^{e,*}$ can be further reduced to enforce
the following localized version of the entropy stability
condition employed in \cite{chen,entropyCG,tadmor}.

\begin{thm}[Entropy correction via subcell flux limiting]\label{thm3}
  Let $\{\eta(u),\mathbf{q}(u)\}$ be an entropy pair
  and $v(u)=\eta'(u)$ the corresponding
  entropy variable. Define $v_h=\sum_{j=1}^{N_h}v_j\varphi_j$ using the
  Bernstein coefficients $v_j=v(u_j),\ j=1,\ldots,N_h$.
  Suppose that the limited subcell fluxes
  $\bar f_{ij}^e$ satisfy 
\beq\label{condES}
\frac{v_i-v_j}{2}
     [\tilde d_{ij}^e(u_j-u_i)+\bar f_{ij}^e
       -\tilde{\mathbf{c}}_{ij}^e\cdot(\mathbf{f}(u_j)+\mathbf{f}(u_i))]
   -\tilde{\mathbf{c}}_{ij}^e\cdot[\boldsymbol{\psi}(u_j)-\boldsymbol{\psi}(u_i)]
    \le \tilde p_{ij}^e,
    \eeq
    where $\boldsymbol{\psi}(u)=v(u)\mathbf{f}(u)-\mathbf{q}(u)$ and
    \beq\label{condES2}
  \sum_{i\in\mathcal N^e}\sum_{j\in\tilde{\mathcal N}^e_i\backslash\{i\}}\tilde p_{ij}^e\le    
  \sum_{i\in\mathcal N^e}\sum_{j\in\mathcal N^e\backslash\{i\}}
  (\tilde{\mathbf{c}}_{ij}^e-\mathbf{c}_{ij}^e)\cdot
  \left[
    \left(\frac{v_i-v_j}{2}\right) (\mathbf{f}(u_j)-\mathbf{f}(u_i)) + \bfq(u_j)-\bfq(u_i)  
    \right]
  =:p^{e,\max}.
  \eeq
Then the flux-corrected semi-discrete scheme \eqref{fcorr}
satisfies the discrete entropy inequality
\beq
\sum_{e=1}^{E_h}\sum_{i\in\mathcal N^e}m_{i}^e\td{\eta(u_i)}{t}\le
\sum_{e=1}^{E_h}\sum_{i\in\mathcal N^e}
\sum_{j\in\mathcal N^e\backslash\{i\}}
   [G_{ij}^e-
      \mathbf{c}_{ij}^e\cdot(\mathbf{q}(u_j)-\mathbf{q}(u_i))],
\eeq
where 
\beq
G_{ij}^e=\frac{v_i+v_j}{2}[\tilde d_{ij}^e(u_j-u_i)+\bar f_{ij}^e]-
  \frac{v_i-v_j}{2}\,\mathbf{c}_{ij}^e\cdot (\mathbf{f}(u_j)-\mathbf{f}(u_i)).
\eeq
  
\end{thm}

\begin{pf} Using the chain rule to differentiate $\eta(u)$
  and substituting
  $v_i$ for $\eta'(u_i)$,
  we obtain the identity  $$\sum_{e\in\mathcal E_i}m_{i}^e\td{\eta(u_i)}{t}
  =\sum_{e\in\mathcal E_i}m_{i}^ev_i\td{u_i}{t}
  =\sum_{e\in\mathcal E_i}v_i\sum_{j\in\tilde{\mathcal N}_i^e\backslash\{i\}}[\tilde g_{ij}^{e}
    -\tilde{\mathbf{c}}_{ij}^e\cdot(\mathbf{f}(u_j)-\mathbf{f}(u_i))],$$
  where $$\tilde g_{ij}^{e}=\tilde d_{ij}^e(u_j-u_i)+\bar f_{ij}^e.$$ By
  definition, the entries $\tilde{\mathbf{c}}_{ij}^e$ of the element
  matrix $\tilde{\mathbf{C}}^e=M_L^e(M_C^e)^{-1}\mathbf{C}^e$
  satisfy the zero sum condition
  $\sum_{j\in\tilde{\mathcal N}_i^e}\tilde{\mathbf{c}}_{ij}^e=\mathbf{0}.$
Following the proof of
   Theorem 1 in \cite{entropyCG}, we use this zero sum property and the
  entropy stability condition \eqref{condES} to estimate the
  rate of entropy production in element $K^e$ as follows:
   \begin{align*}
&  v_i\sum_{j\in\tilde{\mathcal N}_i^e\backslash\{i\}}[\tilde g_{ij}^{e}
    -\tilde{\mathbf{c}}_{ij}^e\cdot(\mathbf{f}(u_j)-\mathbf{f}(u_i))]
  = v_i\sum_{j\in\tilde{\mathcal N}_i^e\backslash\{i\}}[\tilde g_{ij}^{e}-\tilde{\mathbf{c}}_{ij}^e\cdot(\mathbf{f}(u_j)+\mathbf{f}(u_i))]-2v_i\tilde{\mathbf{c}}_{ii}^e\cdot\mathbf{f}(u_i)\\
  &= \sum_{j\in\tilde{\mathcal N}_i^e\backslash\{i\}}\left(
\frac{v_i+v_j}{2}
     [\tilde g_{ij}^{e}-\tilde{\mathbf{c}}_{ij}^e\cdot(\mathbf{f}(u_j)+\mathbf{f}(u_i))]+\frac{v_i-v_j}{2}
     [\tilde g_{ij}^{e}-\tilde{\mathbf{c}}_{ij}^e\cdot(\mathbf{f}(u_j)+\mathbf{f}(u_i))]\right)
  -2v_i\tilde{\mathbf{c}}_{ii}^e\cdot\mathbf{f}(u_i)\\
  &\le \sum_{j\in\tilde{\mathcal N}_i^e\backslash\{i\}}\left(\frac{v_i+v_j}{2}[\tilde g_{ij}^{e}
    -\tilde{\mathbf{c}}_{ij}^e\cdot(\mathbf{f}(u_j)+\mathbf{f}(u_i))]
+\tilde{\mathbf{c}}_{ij}^e\cdot
[\boldsymbol{\psi}(u_j)-\boldsymbol{\psi}(u_i)]+\tilde p_{ij}^e
\right)-2v_i\tilde{\mathbf{c}}_{ii}^e\cdot\mathbf{f}(u_i)\\
&=\sum_{j\in\tilde{\mathcal N}_i^e\backslash\{i\}}\left(\frac{v_i+v_j}{2}[\tilde g_{ij}^{e}
  -\tilde{\mathbf{c}}_{ij}^e\cdot(\mathbf{f}(u_j)+\mathbf{f}(u_i))]
+\tilde{\mathbf{c}}_{ij}^e\cdot
[\boldsymbol{\psi}(u_j)+\boldsymbol{\psi}(u_i)]+\tilde p_{ij}^e\right) -2\tilde{\mathbf{c}}_{ii}^e\cdot[
  v_i\mathbf{f}(u_i)-\boldsymbol{\psi}(u_i)]\\
&=\sum_{j\in\tilde{\mathcal N}_i^e\backslash\{i\}}\left(
\frac{v_i+v_j}{2}\tilde g_{ij}^{e}
- \tilde{\mathbf{c}}_{ij}^e\cdot\left[
  \frac{v_i-v_j}{2}(\mathbf{f}(u_j)-\mathbf{f}(u_i))+
 \mathbf{q}(u_j)+\mathbf{q}(u_i)\right]+\tilde p_{ij}^e\right)
-2\tilde{\mathbf{c}}_{ii}^e\cdot\mathbf{q}(u_i)\\
&=\sum_{j\in\tilde{\mathcal N}_i^e\backslash\{i\}}\left(
\frac{v_i+v_j}{2}\tilde g_{ij}^{e}- \tilde{\mathbf{c}}_{ij}^e\cdot\left[
  \frac{v_i-v_j}{2}(\mathbf{f}(u_j)-\mathbf{f}(u_i))+
  \mathbf{q}(u_j)-\mathbf{q}(u_i)\right]+\tilde p_{ij}^e\right).
\end{align*}
   Summing over $e=1,\ldots,E_h$ and $i\in\mathcal N^e$, we use
   assumption \eqref{condES2} to eliminate the auxiliary quantities
   $\tilde p_{ij}^e$ and obtain
   $\mathbf{c}_{ij}^e$ instead of
   $\tilde{\mathbf{c}}_{ij}^e$ on the right-hand side of the
   final estimate
   \begin{align*}
\sum_{e=1}^{E_h}\sum_{i\in\mathcal N^e}m_{i}^e\td{\eta(u_i)}{t}
&\le   \sum_{e=1}^{E_h}  \sum_{i\in\mathcal N^e}
\sum_{j\in\mathcal N^e\backslash\{i\}}\left(
\frac{v_i+v_j}{2}\tilde g_{ij}^{e}-\mathbf{c}_{ij}^e\cdot\left[
  \frac{v_i-v_j}{2}(\mathbf{f}(u_j)-\mathbf{f}(u_i))+
  \mathbf{q}(u_j)-\mathbf{q}(u_i)\right]\right)\\&=
\sum_{e=1}^{E_h}\sum_{i\in\mathcal N^e}
\sum_{j\in\mathcal N^e\backslash\{i\}}
   [G_{ij}^e-
      \mathbf{c}_{ij}^e\cdot(\mathbf{q}(u_j)-\mathbf{q}(u_i))]
   \end{align*}
  which proves the assertion of the theorem.
\qquad\proofend

\end{pf}

\begin{rmk}
  In view of property \eqref{skew},
  a further rearrangement yields
  the estimate (cf. \cite{entropyCG})
\begin{align}
\sum_{e=1}^{E_h}\sum_{i\in\mathcal N^e}
 m_i^e\td{\eta(u_i)}{t}&\le\sum_{i=1}^{N_h}\Big(\sum_{j\in\mathcal N_i\atop j>i}
 \frac{v_i+v_j}{2}\sum_{e\in\mathcal E_i\cap\mathcal E_j}\underbrace{(\tilde
   g_{ij}^{e}+\tilde g_{ji}^{e})}_{=0}-2\mathbf{q}(u_i)
\cdot\underbrace{\sum_{e\in\mathcal E_i}\mathbf{c}_{ii}^e}_{=\mathbf{0}}\Big)
  \nonumber \\
&-\sum_{i=1}^{N_h}
\sum_{j\in\mathcal N_i\atop j>i}
\left[\frac{v_i-v_j}{2}(\mathbf{f}(u_j)-\mathbf{f}(u_i))+
  \mathbf{q}(u_j)+\mathbf{q}(u_i)\right]\cdot
\underbrace{\sum_{e\in\mathcal E_i\cap\mathcal E_j}(\mathbf{c}_{ij}^e+\mathbf{c}_{ji}^e)}_{=0}
=0\label{detadt}
\end{align}
in accordance with the fact that $\td{}{t}\int_{\Omega}\eta(u)\dx
\le 0$ under the assumption of periodic boundary conditions.
Note that condition \eqref{skew} and inequality \eqref{detadt}
do not hold if the coefficients $\mathbf{c}_{ij}^e$ are
replaced with the coefficients $\tilde{\mathbf{c}}_{ij}^e$ of
the lumped discrete gradient. That is why we impose condition
\eqref{condES2} and use it to obtain the final estimate
in terms of $\mathbf{c}_{ij}^e$ rather than $\tilde{\mathbf{c}}_{ij}^e$
in the proof of Theorem \ref{thm3}.
  \end{rmk}

For practical limiting purposes, we still need to define (i) distributed production
bounds $\tilde p_{ij}^e$ such that assumption \eqref{condES2} holds and (ii)
subcell fluxes $\bar f_{ij}^e$ such that the entropy stability condition
\eqref{condES} holds. The low-order LLF approximation corresponding
to \eqref{fcorr} with $\bar f_{ij}^e=0$ satisfies \cite{chen,entropyCG}
\beq\label{q_tilde}
\tilde q_{ij}^e:=\frac{v_i-v_j}{2}
     [\tilde d_{ij}^e(u_j-u_i)
     -\tilde{\mathbf{c}}_{ij}^e\cdot(\mathbf{f}(u_j)+\mathbf{f}(u_i))]
   -\tilde{\mathbf{c}}_{ij}^e\cdot[\boldsymbol{\psi}(u_j)-\boldsymbol{\psi}(u_i)]
   \le 0.
   \eeq
To ensure that condition \eqref{condES} holds for $\bar f_{ij}^e=0$ and,
therefore, can be enforced by adjusting the magnitude of $\bar f_{ij}^e$,
we must have $\tilde q_{ij}^e \le \tilde p_{ij}^e$. 
If $p^{e,\max}\ge 0$, we set $\tilde p_{ij}^e=0$ for all
$j\in\tilde{\mathcal N_i}^e\backslash\{i\}$. A negative
 entropy production bound
 $p^{e,\max}$  can be split into a sum of components $\tilde p_{ij}^e$
 as follows:
 \beq\label{p_tilde}
\tilde p_{ij}^e=\omega_{ij}^e\min\{0,p^{e,\max}\},\qquad
\omega_{ij}^e=\frac{\tilde q_{ij}^e-\epsilon}{\sum_{k\in\mathcal N^e}\sum_{l\in\tilde{\mathcal N}^e_k\backslash\{k\}}(\tilde q_{kl}^e-\epsilon)},
\eeq
where $\epsilon>0$ is an infinitesimally small positive number that
we use to formally prevent division by zero. After distributing 
$p^{e,\max}$ among pairs of nearest neighbor nodes
in this way, we check the validity of the feasibility
condition $\tilde q_{ij}^e \le \tilde p_{ij}^e$
which can always be enforced by adding
\beq\label{dAdd}
  \tilde d_{ij}^{e,{\rm add}}=\frac{2\min
    \{\tilde p_{ij}^e-\tilde q_{ij}^e, 0, \tilde p_{ji}^e-\tilde q_{ji}^e\}
  }{(v_i-v_j)(u_j-u_i)-\epsilon}
\eeq
to $\tilde d_{ij}^{e}$ if necessary. By the mean value theorem,
we have $v_i-v_j=\eta'(u_i)-\eta'(u_j)
=\eta''(\xi)(u_i-u_j)$  for some $\xi\in\R$. It follows that
$(v_i-v_j)(u_j-u_i)\le 0$ and, therefore, $\tilde d_{ij}^{e,\rm add}\ge 0$
for any convex entropy $\eta$.
\medskip

The bound-preserving flux $\tilde f_{ij}^{e,*}$, as
defined by \eqref{fij_lim}, can now be
adjusted to satisfy \eqref{condES} as follows:
\beq\label{fij_lim_fix}
\bar f_{ij}^{e}
=\begin{cases}
\frac{\min\{2\bar p_{ij}^{e},(v_i-v_j)\tilde f_{ij}^{e,*},2\bar p_{ji}^{e}\}}{v_i-v_j}
&\mbox{if}\ (v_i-v_j)\tilde f_{ij}^{e,*}>0,\\
 \tilde f_{ij}^* & \mbox{otherwise},
\end{cases}
\eeq
where
\beq\label{bar_p}
\bar p_{ij}^{e}=\tilde p_{ij}^e-\tilde q_{ij}^{e}-\frac{v_i-v_j}{2}\tilde d_{ij}^{e,\rm add}(u_j-u_i)\ge 0
\eeq
are nonnegative
upper bounds for entropy-producing subcell fluxes. Since the limiting
procedure is similar to that developed in \cite{entropyCG} for
$\mathbb{P}_1/\mathbb{Q}_1$ elements, we refer the reader to \cite{entropyCG}
for further details.

\section{Summary of the algorithm}\label{sec:summary}
For the reader's convenience, we summarize the proposed algorithm in this
section. If no flux limiters are applied, then the method is given by \eqref{CGstabEV}. Otherwise, we proceed as follows: 
\begin{itemize}
\item[1.] Compute the Bernstein coefficients $\dot{u}^S_i, ~i=1, \ldots, N_h$ of the stabilized time derivatives via \eqref{targetudot}. 
\item[2.] For each element $K^e,\ e=1,\ldots,E_h$, assemble and solve the subcell system \eqref{aux} 
  to obtain the subcell flux potentials $\dot{w}^e_i, ~i\in\mathcal N^e$. 
\item[3.] Calculate the subcell fluxes 
 $\tilde f_{ij}^e$ defined by \eqref{subcell} and their IDP
 counterparts $\tilde f_{ij}^{e,*}$  defined by
 \eqref{fij_lim}. 
\end{itemize}

If no additional entropy check is performed, then the flux-corrected
semi-discrete scheme is given by 
\eqref{fcorr} with $\bar f_{ij}^e=\tilde f_{ij}^{e,*}$. Otherwise, the
following steps complete the process of subcell flux limiting:
\begin{itemize}
\item[4.] Compute $\tilde q_{ij}^e$ and $\tilde p_{ij}^e$ via \eqref{q_tilde} and 
  \eqref{p_tilde}, respectively.
\item[5.] 
  If $\tilde q_{ij}^e>\tilde p_{ij}^e$, calculate
  $\tilde d_{ij}^{e,\text{add}}$ via \eqref{dAdd}. Otherwise, set
   $\tilde d_{ij}^{e,\text{add}}:=0$.
\item[6.] Calculate the entropy stability preserving flux
  $\bar f_{ij}^{e,*}$ via \eqref{fij_lim_fix},\eqref{bar_p}.

\end{itemize}
The entropy-corrected semi-discrete scheme is then given by \eqref{fcorr} with
$\tilde d_{ij}^e$ replaced by $\tilde d_{ij}^e+\tilde d_{ij}^{e,\text{add}}$. 

\begin{rmk}\label{remark:linear_adv}
  For linear advection with constant velocity $\mathbf{v}\in\R^d$, we have
  $\mathbf{f}(u)=\mathbf{v}u$. In this case, the entropy stability condition
 based on  $\eta(u)=\frac{u^2}2$ and  $\mathbf{q}(u)=\frac12\mathbf{v}u^2$
  reduces to
\beq
\frac{u_i-u_j}{2}
     [(\tilde d_{ij}^e+\tilde d_{ij}^{e,\rm add})(u_j-u_i)+\bar f_{ij}^e]
    \le \tilde p_{ij}^e\le 0,
    \eeq
    which means that the net flux $(\tilde d_{ij}^e+\tilde d_{ij}^{e,\rm add})(u_j-u_i)
    +\bar f_{ij}^e$ must be entropy-dissipative to satisfy \eqref{condES} with
    $\tilde d_{ij}^e+\tilde d_{ij}^{e,\rm add}$ in place of $\tilde d_{ij}^e$. As
    we show in Section \ref{sec:num}, entropy limiting based on this criterion
    may increase the levels of numerical dissipation without having any positive
    effect in the case of linear advection. Hence, it is worthwhile to omit
    steps 4-6 in applications to linear advection problems.
    
\end{rmk}  

\section{Numerical examples}\label{sec:num}
In this section, we perform numerical experiments for linear and nonlinear
scalar problems. Our numerical examples illustrate the impact of each
correction step (linear stabilization, nonlinear stabilization, IDP
limiting, entropy fix) and the properties of resulting methods in
different situations. For test problems with smooth solutions, we
show that optimal convergence behavior can be achieved with
stabilized high-order schemes presented in Sections \ref{sec:ls_stab}
and \ref{sec:ev_stab}, respectively. The results for nonlinear problems
with shocks and nonconvex flux functions demonstrate the ability of
the limiting procedure proposed in Section \ref{sec:ES-lim}
to prevent spurious oscillations and convergence to wrong weak
solutions.

In the description of our numerical results, the methods under
investigation are labeled as follows:
\begin{itemize}
\item HO-X: high-order Galerkin method equipped with linear
  stabilization of type $\mathrm{X}\in\{{\rm SUPG, VMS}\}$
  (as defined in Section \ref{sec:ls_stab}: 
  no nonlinear entropy stabilization, no convex limiting);
\item HO-X-EV: entropy-stabilized counterpart of HO-X 
 (as defined in Section \ref{sec:ev_stab}: no convex limiting);
\item Y-BP: bound-preserving counterpart of  
  $\text{Y}\in\{\text{HO-X}, ~\text{HO-X-EV}\}$ corresponding
 to the algebraic flux correction scheme \eqref{fcorr} 
  with 
  $\bar f_{ij}^e=\tilde{f}_{ij}^{e,*}$, where $\tilde{f}_{ij}^{e,*}$ 
  is given by \eqref{fij_lim};
\item Y-FL: flux-limited counterpart \eqref{fcorr} of  $\text{Y}\in\{\text{HO-X}, ~\text{HO-X-EV}\}$ with $\bar f_{ij}^e$ defined by \eqref{fij_lim_fix}.
\end{itemize}
In all nonlinear stabilization terms, we use the square entropy $\eta(u)=\frac{u^2}2$. 
Discretization in time is performed using the third-order
explicit SSP Runge-Kutta method with three stages \cite{ssprev}
unless mentioned otherwise.
The implementation of all methods
is based on the open-source C++  library MFEM \cite{mfem}
which provides optimized tools for computations with
high-order finite elements.

\subsection{Linear advection with constant velocity in 1D}

To determine the experimental order of convergence (EOC) for the
stabilized high-order methods HO-X and HO-X-EV, we apply them to
the one-dimensional linear advection equation 
\beq\label{linear_advection}
\pd{u}{t}+a\pd{u}{x}=0\quad\mbox{in}\quad\Omega=(0,1)
\eeq
with constant velocity $a=1$. The first initial condition that we
consider is given by
\beq\label{initcos}
u_0(x) = \cos(2\pi(x-0.5)).
\eeq
We evolve this smooth profile up to a final time $t=1$ on a sequence of
successively refined uniform grids and measure the EOCs w.r.t. the $L^1$ norm. 
To keep the temporal errors negligible, we discretize in time using
a 6th order explicit
Runge-Kutta method whose Butcher tableau
 is given by \cite{rkmethods}
\begin{align}\label{RK76}
\begin{tabular}{c|ccccccc}
  0 & & & & & & & \\ 
  1/3 & 1/3 & & & & & & \\
  2/3 & 0 & 2/3 & & & & & \\
  1/3 & 1/12 & 1/3 & -1/12 & & & & \\
  1/2 & -1/16 & 9/8 & -3/16 & -3/8 & & & \\
  1/2 & 0 & 9/8 & -3/8 & -3/4 & 1/2 & & \\
  1   & 9/44 & -9/11 & 63/44 & 18/11 & 0 & -16/11 \\ \hline
  & 11/120 & 0 & 27/40 & 27/40 & -4/15 & -4/15 & 11/120
\end{tabular}
\end{align}
The results of the grid convergence study
are shown in Table \ref{table:advection}. 
All methods deliver the optimal EOCs. 

\begin{table}[!h]\scriptsize
  \begin{center}
    \subfloat[$p=1$]{
    \begin{tabular}{||c||c|c||c|c||c|c||c|c||} \cline{1-9}
      \multicolumn{1}{||c||}{} &
      \multicolumn{2}{|c||}{\bf HO-SUPG} &
      \multicolumn{2}{|c||}{\bf HO-VMS} &
      \multicolumn{2}{|c||}{\bf HO-SUPG-EV} &
      \multicolumn{2}{|c||}{\bf HO-VMS-EV} \\ \cline{2-9}
      $N_h$ &
      $\|u_h-u_\text{exact}\|_{L^1}$ & EOC &
      $\|u_h-u_\text{exact}\|_{L^1}$ & EOC  &
      $\|u_h-u_\text{exact}\|_{L^1}$ & EOC  &
      $\|u_h-u_\text{exact}\|_{L^1}$ & EOC \\ \hline
    8   & 8.46E-02 &  --  & 2.02E-01 &  --  & 8.46E-02 &  --  & 2.02E-01 &  --  \\ 
    16  & 1.03E-02 & 3.03 & 2.96E-02 & 2.77 & 1.03E-02 & 3.03 & 2.96E-02 & 2.77 \\
    32  & 1.27E-03 & 3.02 & 3.78E-03 & 2.97 & 1.27E-03 & 3.02 & 3.78E-03 & 2.97 \\
    64  & 2.10E-04 & 2.60 & 4.73E-04 & 3.00 & 2.10E-04 & 2.60 & 4.73E-04 & 3.00 \\
    128 & 5.46E-05 & 1.94 & 5.92E-05 & 3.00 & 5.46E-05 & 1.94 & 5.92E-05 & 3.00 \\
    256 & 1.39E-05 & 1.97 & 1.34E-05 & 2.14 & 1.39E-05 & 1.97 & 1.34E-05 & 2.14 \\ \hline
    \end{tabular}
      }

    \subfloat[$p=2$]{
    \begin{tabular}{||c||c|c||c|c||c|c||c|c||} \cline{1-9}
      \multicolumn{1}{||c||}{} &
      \multicolumn{2}{|c||}{\bf HO-SUPG} &
      \multicolumn{2}{|c||}{\bf HO-VMS} &
      \multicolumn{2}{|c||}{\bf HO-SUPG-EV} &
      \multicolumn{2}{|c||}{\bf HO-VMS-EV} \\ \cline{2-9}
      $N_h$ &
      $\|u_h-u_\text{exact}\|_{L^1}$ & EOC &
      $\|u_h-u_\text{exact}\|_{L^1}$ & EOC  &
      $\|u_h-u_\text{exact}\|_{L^1}$ & EOC  &
      $\|u_h-u_\text{exact}\|_{L^1}$ & EOC \\ \hline
      16  & 1.95E-03 &  --  & 1.98E-03 &  --  & 8.03E-03 &  --  & 7.52E-03 &  --  \\
      32  & 2.28E-04 & 3.09 & 2.28E-04 & 3.12 & 5.36E-04 & 3.91 & 5.24E-04 & 3.84 \\
      64  & 2.77E-05 & 3.04 & 2.77E-05 & 3.04 & 4.63E-05 & 3.53 & 4.86E-05 & 3.43 \\
      128 & 3.44E-06 & 3.01 & 3.44E-06 & 3.01 & 4.91E-06 & 3.23 & 5.07E-06 & 3.26 \\
      256 & 4.29E-07 & 3.00 & 4.29E-07 & 3.00 & 5.58E-07 & 3.14 & 5.74E-07 & 3.14 \\
      512 & 5.36E-08 & 3.00 & 5.36E-08 & 3.00 & 6.53E-08 & 3.09 & 6.91E-08 & 3.05 \\ \hline
    \end{tabular}
      }

    \subfloat[$p=3$]{
    \begin{tabular}{||c||c|c||c|c||c|c||c|c||} \cline{1-9}
      \multicolumn{1}{||c||}{} &
      \multicolumn{2}{|c||}{\bf HO-SUPG} &
      \multicolumn{2}{|c||}{\bf HO-VMS} &
      \multicolumn{2}{|c||}{\bf HO-SUPG-EV} &
      \multicolumn{2}{|c||}{\bf HO-VMS-EV} \\ \cline{2-9}
      $N_h$ &
      $\|u_h-u_\text{exact}\|_{L^1}$ & EOC &
      $\|u_h-u_\text{exact}\|_{L^1}$ & EOC  &
      $\|u_h-u_\text{exact}\|_{L^1}$ & EOC  &
      $\|u_h-u_\text{exact}\|_{L^1}$ & EOC \\ \hline
      24  & 1.14E-04 &  --  & 1.52E-04 &  --  & 1.85E-04 &  --  & 2.82E-04 &  --  \\ 
      48  & 6.36E-06 & 4.16 & 9.77E-06 & 3.96 & 8.02E-06 & 4.52 & 1.35E-05 & 4.38 \\
      96  & 3.79E-07 & 4.07 & 6.14E-07 & 3.99 & 4.27E-07 & 4.23 & 7.19E-07 & 4.23 \\
      192 & 2.34E-08 & 4.02 & 3.88E-08 & 3.98 & 2.48E-08 & 4.11 & 4.18E-08 & 4.11 \\
      384 & 1.45E-09 & 4.01 & 2.45E-09 & 3.99 & 1.49E-09 & 4.05 & 2.53E-09 & 4.05 \\
      768 & 9.08E-11 & 4.00 & 1.54E-10 & 3.99 & 9.18E-11 & 4.02 & 1.56E-10 & 4.02 \\ \hline
    \end{tabular}
      }

    \subfloat[$p=4$]{
    \begin{tabular}{||c||c|c||c|c||c|c||c|c||} \cline{1-9}
      \multicolumn{1}{||c||}{} &
      \multicolumn{2}{|c||}{\bf HO-SUPG} &
      \multicolumn{2}{|c||}{\bf HO-VMS} &
      \multicolumn{2}{|c||}{\bf HO-SUPG-EV} &
      \multicolumn{2}{|c||}{\bf HO-VMS-EV} \\ \cline{2-9}
      $N_h$ &
      $\|u_h-u_\text{exact}\|_{L^1}$ & EOC &
      $\|u_h-u_\text{exact}\|_{L^1}$ & EOC  &
      $\|u_h-u_\text{exact}\|_{L^1}$ & EOC  &
      $\|u_h-u_\text{exact}\|_{L^1}$ & EOC \\ \hline
      32   & 2.94E-06 &  --  & 3.09E-06 &  --  & 3.84E-06 &  --  & 5.24E-06 &  --  \\
      64   & 8.62E-08 & 5.09 & 8.98E-08 & 5.11 & 1.03E-07 & 5.22 & 1.21E-07 & 5.44 \\ 
      128  & 2.63E-09 & 5.03 & 2.77E-09 & 5.02 & 2.95E-09 & 5.12 & 3.28E-09 & 5.20 \\
      256  & 8.20E-11 & 5.00 & 8.66E-11 & 5.00 & 8.77E-11 & 5.07 & 9.54E-11 & 5.10 \\
      512  & 2.57E-12 & 5.00 & 2.69E-12 & 5.01 & 2.66E-12 & 5.04 & 2.84E-12 & 5.07 \\ \hline
    \end{tabular}
      }
    \caption{One-dimensional
      linear advection problem \eqref{linear_advection} with
  the initial condition \eqref{initcos}. 
  Grid convergence history for stabilized
  high-order methods using finite elements of degree $p\in\{1,2,3,4\}$.
  \label{table:advection}}
  \end{center} 
\end{table}

Let us now examine
the long time behavior of the high-order entropy-stabilized HO-X-EV
methods for the linear advection problem \eqref{linear_advection}
 with initial data given by \cite{Guermond2011}
\beq\label{ltime_behavior}
u_0(x) = 
\begin{cases}
  e^{-300(2x-0.3)^2} & \mbox{ if } |2x-0.3|\leq 0.25, \\
  1                  & \mbox{ if } |2x-0.9|\leq 0.2, \\
  \sqrt{1-\left(\frac{2x-1.6}{0.2}\right)^2} & \mbox{ if } |2x-1.6|\leq 0.2, \\
  0 & \mbox{\text{otherwise}}.
\end{cases}
\eeq

\begin{figure}[!h]
  \centering
  \includegraphics[scale=0.5]{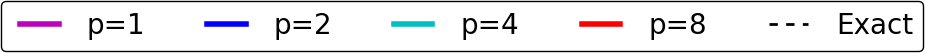}

  \subfloat[standard CG method \eqref{CGhigh}\label{fig:long_time_behavior_pure_galerkin}]{
    \includegraphics[scale=0.325]{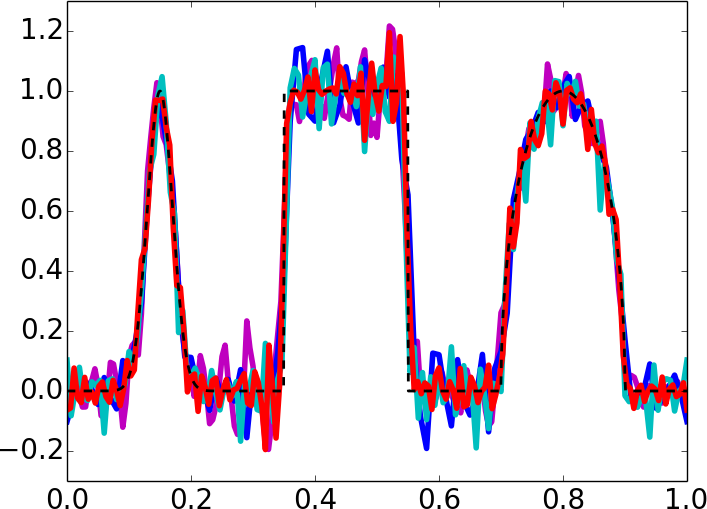}}
  \qquad
  \subfloat[HO-VMS with $\omega=0.1$\label{fig:long_time_behavior_ho_vms_t100}]{
    \includegraphics[scale=0.325]{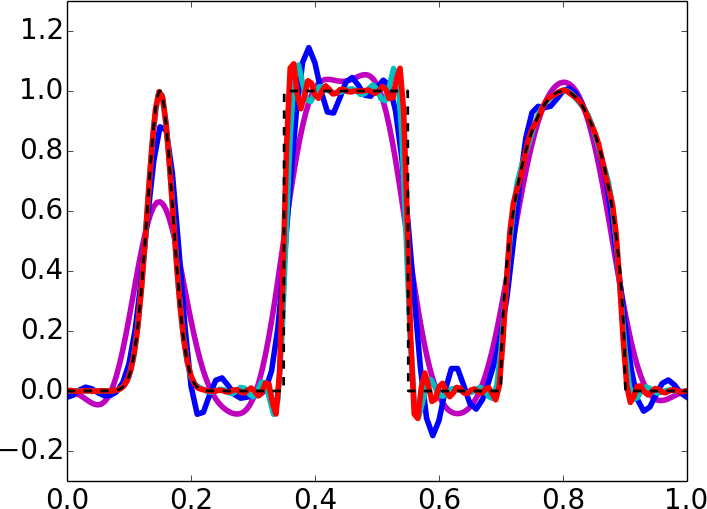}}

  \subfloat[HO-VMS-EV with $\omega=0$ \label{fig:long_time_behavior_pure_ent_visc_t100}]{
    \includegraphics[scale=0.325]{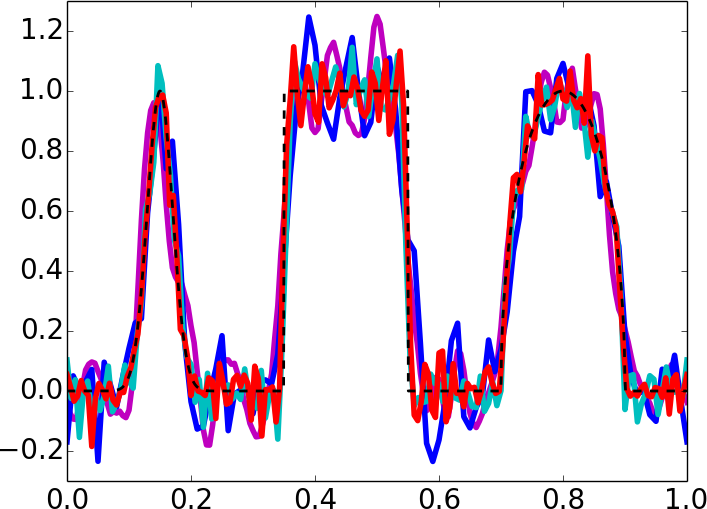}}
  \qquad
  \subfloat[HO-VMS-EV with $\omega=0.1$\label{fig:long_time_behavior_both_t100}]{
    \includegraphics[scale=0.325]{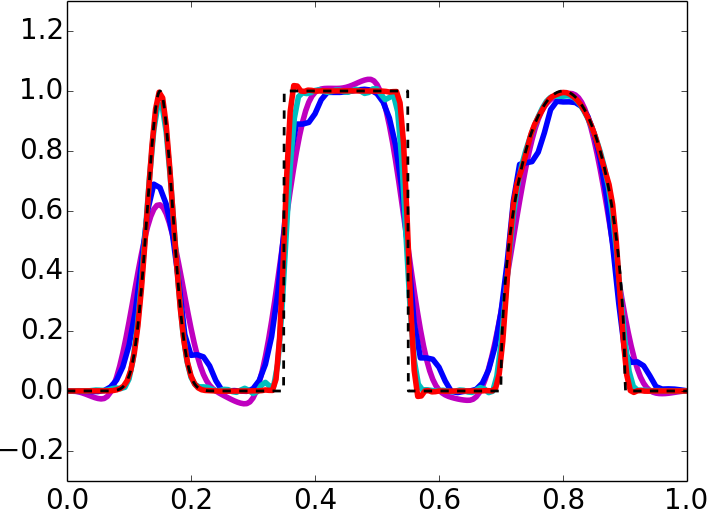}}
  \caption{One-dimensional linear advection problem \eqref{linear_advection}
    with initial condition  \eqref{ltime_behavior}. 
    Numerical solutions at $t=100$ obtained with stabilized high-order methods 
    using $N_h=200$ DoFs.}
  \label{fig:long_time_behavior_t100}
\end{figure}

We run the simulations up to the final time $t=100$ using piecewise-polynomial 
finite element spaces of degree $p=\{1,2,4,8\}$. The mesh size is chosen in such a way 
that the total number of degrees of freedom (DoFs) is $N_h=200$ for each space. 
That is, our high-order finite element approximations use larger mesh cells 
than their low-order counterparts. 
In Figure \ref{fig:long_time_behavior_pure_galerkin}, we show the numerical solutions obtained with the standard Galerkin method \eqref{CGhigh}. As expected, these solutions are highly oscillatory. The amplitude of spurious oscillations can be greatly reduced by using any of the high-order linear stabilization techniques presented in Section \ref{sec:ls_stab}. In Figure \ref{fig:long_time_behavior_ho_vms_t100}, we present the results produced by HO-VMS with $\omega=0.1$. The numerical solutions obtained without any linear stabilization (i.e., using $\omega=0$) are shown in Fig.~\ref{fig:long_time_behavior_pure_ent_visc_t100}. It can be seen that EV stabilization alone is insufficient for linear problems. On
the other hand, Figs \ref{fig:long_time_behavior_ho_vms_t100}
and \ref{fig:long_time_behavior_both_t100} demonstrate that
HO-VMS-EV exhibits better discontinuity-capturing properties than HO-VMS. The entropy-stabilized numerical solutions are essentially nonoscillatory in this example. The results obtained with HO-SUPG and HO-SUPG-EV are similar (not shown here). The findings of Guermond et al. \cite{Guermond2011} also indicate that methods equipped with nonlinear EV stabilization tend to produce smaller undershoots/overshoots in proximity to steep gradients.

%

\subsection{One-dimensional inviscid Burgers equation}
To study the shock-capturing capabilities of entropy stabilization in the context of nonlinear
problems, we apply our HO-X-EV methods to the inviscid Burgers equation 
\beq\label{1Dburgers}
\pd{u}{t}+\pd{(u^2/2) }{x}=0\quad\mbox{in}\quad\Omega=(0,1).
\eeq
The smooth initial condition is given by 
\beq\label{1Dburgers_init}
u_0(x) = \sin(2\pi x). 
\eeq
The entropy solution of this initial value problem develops
a shock at the critical time $t_c=\frac{1}{2\pi}$. For $t<t_c$,
the smooth exact solution is defined by the nonlinear equation
\beq
u(x,t) = \sin(2\pi(x-u(x,t)t))
\eeq
which can be derived
using the method of characteristics.

We evolve the HO-X and HO-X-EV approximations up to the final time $t=0.1$ and
measure the EOCs w.r.t. the $L^1$ norm. The discretization in time
is again performed using the 6th order Runge-Kutta method with
the Butcher tableau \eqref{RK76}. The results of our grid
convergence studies are presented in Table~\ref{table:1Dburgers}.
All methods under investigation deliver the optimal rates of convergence. 

\begin{table}[!h]\scriptsize
  \begin{center}
    \subfloat[$p=1$]{
    \begin{tabular}{||c||c|c||c|c||} \cline{1-5}
      \multicolumn{1}{||c||}{} &
      \multicolumn{2}{|c||}{\bf HO-SUPG-EV} &
      \multicolumn{2}{|c||}{\bf HO-VMS-EV} \\ \cline{2-5}
      $N_h$ &
      $\|u_h-u_\text{exact}\|_{L^1}$ & EOC  &
      $\|u_h-u_\text{exact}\|_{L^1}$ & EOC \\ \hline
      8   & 5.06E-02 &  --  & 5.42E-02 & --   \\ 
      16  & 9.27E-03 & 2.45 & 1.03E-02 & 2.39 \\ 
      32  & 1.75E-03 & 2.41 & 2.05E-03 & 2.33 \\
      64  & 3.90E-04 & 2.16 & 4.03E-04 & 2.35 \\
      128 & 9.17E-05 & 2.09 & 9.35E-05 & 2.11 \\ 
      256 & 2.25E-05 & 2.03 & 2.28E-05 & 2.04 \\ \hline
    \end{tabular}
      }
    \subfloat[$p=2$]{
    \begin{tabular}{||c||c|c||c|c||} \cline{1-5}
      \multicolumn{1}{||c||}{} &
      \multicolumn{2}{|c||}{\bf HO-SUPG-EV} &
      \multicolumn{2}{|c||}{\bf HO-VMS-EV} \\ \cline{2-5}
      $N_h$ &
      $\|u_h-u_\text{exact}\|_{L^1}$ & EOC  &
      $\|u_h-u_\text{exact}\|_{L^1}$ & EOC \\ \hline
      16  & 5.48E-03 &  --  & 5.41E-03 &  --  \\
      32  & 5.57E-04 & 3.30 & 6.15E-04 & 3.14 \\
      64  & 1.31E-04 & 2.08 & 1.53E-04 & 2.01 \\
      128 & 1.78E-05 & 2.88 & 1.91E-05 & 3.01 \\ 
      256 & 2.29E-06 & 2.96 & 2.32E-06 & 3.04 \\ 
      512 & 2.86E-07 & 3.00 & 2.87E-07 & 3.02 \\ \hline
    \end{tabular}
      }

    \subfloat[$p=3$]{
    \begin{tabular}{||c||c|c||c|c||} \cline{1-5}
      \multicolumn{1}{||c||}{} &
      \multicolumn{2}{|c||}{\bf HO-SUPG-EV} &
      \multicolumn{2}{|c||}{\bf HO-VMS-EV} \\ \cline{2-5}
      $N_h$ &
      $\|u_h-u_\text{exact}\|_{L^1}$ & EOC  &
      $\|u_h-u_\text{exact}\|_{L^1}$ & EOC \\ \hline
      24  & 7.30E-04 &  --  & 7.42E-04 & --   \\
      48  & 2.16E-04 & 1.76 & 2.34E-04 & 1.66 \\
      96  & 1.43E-05 & 3.91 & 1.86E-05 & 3.66 \\
      192 & 7.96E-07 & 4.17 & 1.30E-06 & 3.83 \\ 
      384 & 4.94E-08 & 4.01 & 8.70E-08 & 3.91 \\
      768 & 3.12E-09 & 3.99 & 5.61E-09 & 3.96 \\ \hline
    \end{tabular}
    }
    \subfloat[$p=4$]{
    \begin{tabular}{||c||c|c||c|c||} \cline{1-5}
      \multicolumn{1}{||c||}{} &
      \multicolumn{2}{|c||}{\bf HO-SUPG-EV} &
      \multicolumn{2}{|c||}{\bf HO-VMS-EV} \\ \cline{2-5}
      $N_h$ &
      $\|u_h-u_\text{exact}\|_{L^1}$ & EOC  &
      $\|u_h-u_\text{exact}\|_{L^1}$ & EOC \\ \hline
      32   & 6.03E-04 &  --  & 6.26E-04 &  --  \\
      64   & 3.72E-05 & 4.02 & 4.14E-05 & 3.92 \\
      128  & 4.66E-07 & 6.32 & 7.11E-07 & 5.86 \\
      256  & 3.54E-08 & 3.72 & 4.47E-08 & 3.99 \\
      512  & 1.07E-09 & 5.05 & 1.16E-09 & 5.26 \\
      1024 & 3.28E-11 & 5.03 & 3.45E-11 & 5.07 \\ \hline
    \end{tabular}
    }
    \caption{One-dimensional inviscid Burgers equation \eqref{1Dburgers} with initial condition \eqref{1Dburgers_init}. Grid convergence history for stabilized high-order methods using finite elements of degree $p\in\{1,2,3,4\}$.
      \label{table:1Dburgers}}
  \end{center}
\end{table}

\begin{figure}[!h]
  \centering
    \includegraphics[scale=0.4]{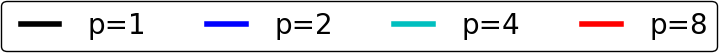}

  \subfloat[HO-SUPG\label{fig:1Dburgers_HO_SUPG}]{
    \includegraphics[scale=0.2]{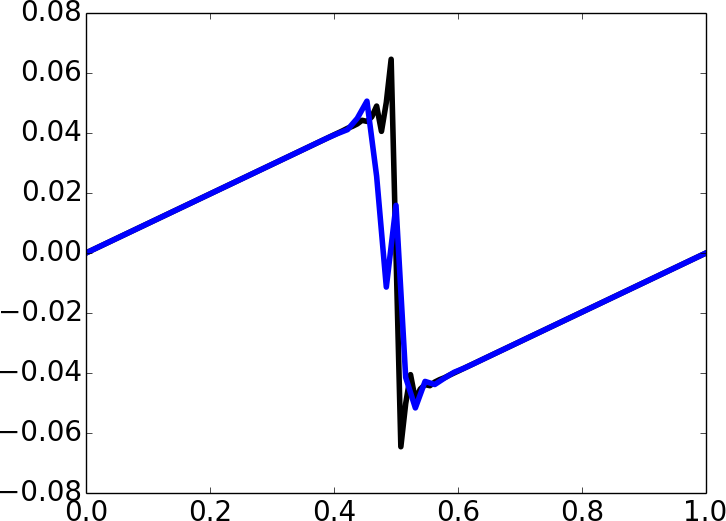}}
  \quad
  \subfloat[HO-VMS\label{fig:1Dburgers_HO_VMS}]{
    \includegraphics[scale=0.2]{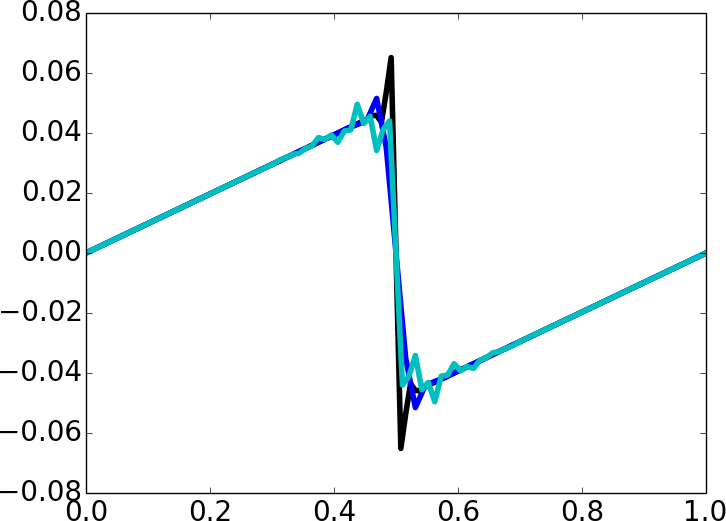}}
  \quad
  \subfloat[HO-SUPG-EV\label{fig:1Dburgers_HO_SUPG_EV}]{
    \includegraphics[scale=0.2]{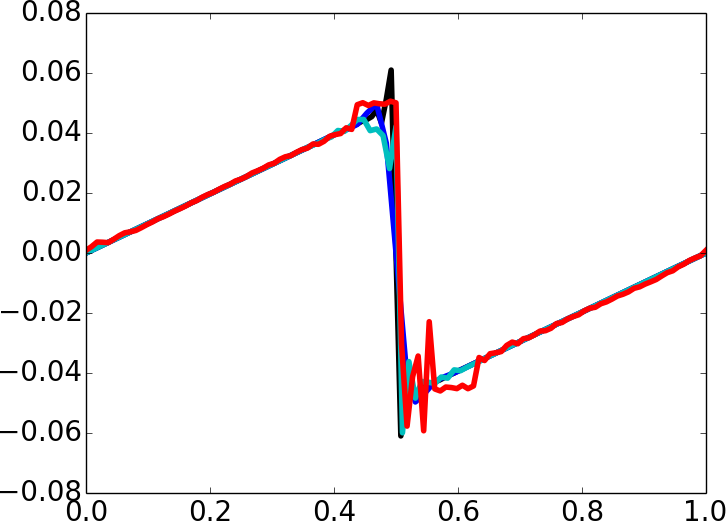}}
  \quad
  \subfloat[HO-VMS-EV\label{fig:1Dburgers_HO_VMS_EV}]{
    \includegraphics[scale=0.2]{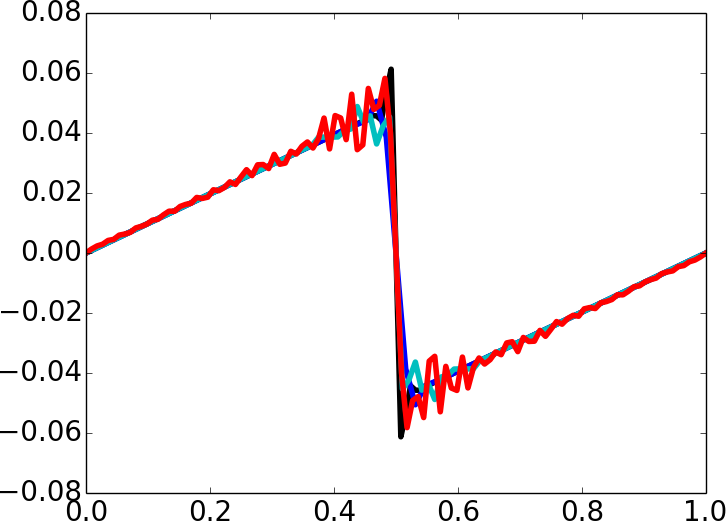}}

  \subfloat[HO-SUPG-FL\label{fig:1Dburgers_MCL_HO_SUPG}]{
    \includegraphics[scale=0.2]{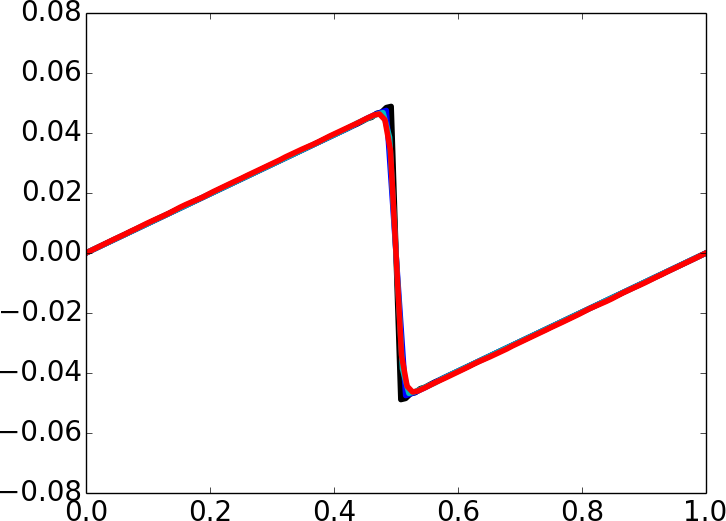}}
  \quad
  \subfloat[HO-VMS-FL\label{fig:1Dburgers_MCL_HO_VMS}]{
    \includegraphics[scale=0.2]{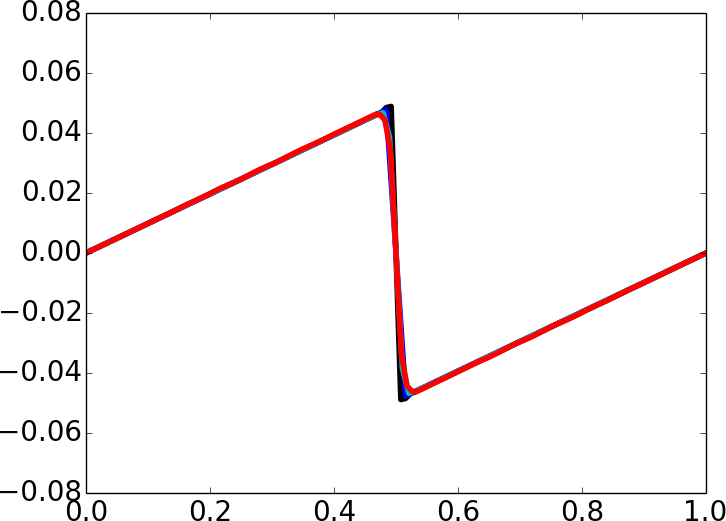}}
  \quad
  \subfloat[HO-SUPG-EV-FL\label{fig:1Dburgers_MCL_HO_SUPG_EV}]{
    \includegraphics[scale=0.2]{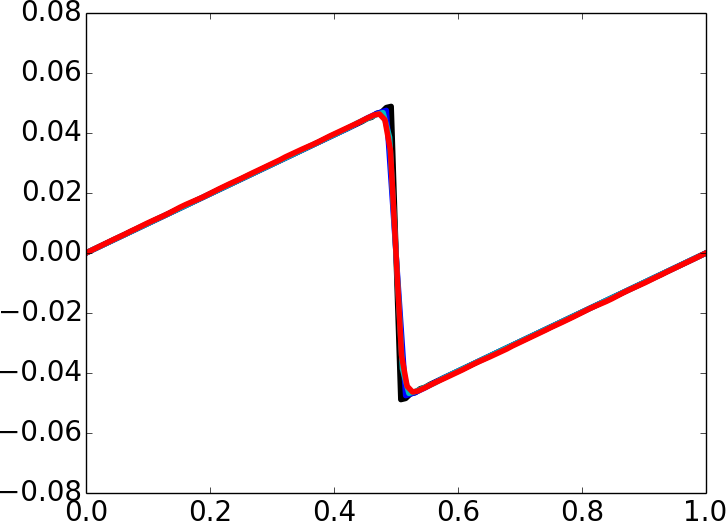}}
  \quad
  \subfloat[HO-VMS-EV-FL\label{fig:1Dburgers_MCL_HO_VMS_EV}]{
    \includegraphics[scale=0.2]{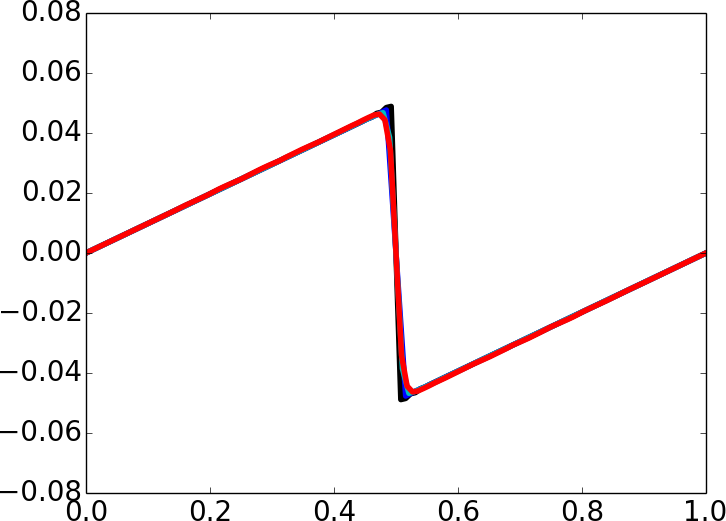}}
  \caption{One-dimensional Burgers equation \eqref{1Dburgers} with initial condition \eqref{1Dburgers_init}. Numerical solutions at $t=10$ obtained using $N_h=128$ DoFs.
  \label{fig:1Dburgers}}
\end{figure}

Let us now run the simulations up to $t=10$ and study the ability of the stabilized HO methods to capture the shock that forms at $t=t_c$. Computations are performed using piecewise-polynomial finite element spaces of degree $p=\{1,2,4,8\}$. In all numerical experiments, we use $N_h=128$ DoFs. The results obtained with HO-X and HO-X-EV are shown in Figs \ref{fig:1Dburgers_HO_SUPG}--\ref{fig:1Dburgers_HO_VMS_EV}. We remark that our simulations became unstable for HO-SUPG with $p=\{4,8\}$ and HO-VMS with $p=8$. The use of entropy stabilization has cured this problem. In addition, it reduced the magnitude of
spurious oscillations in all cases. The nonoscillatory solutions
shown in Figs \ref{fig:1Dburgers_MCL_HO_SUPG}--\ref{fig:1Dburgers_MCL_HO_VMS_EV}
were obtained using the monolithic convex limiting techniques of Section
\ref{sec:ES-lim}. It can be seen that the flux-limited (FL) version of
each high-order method enforces local bounds without introducing
large amounts of numerical diffusion. The differences between the FL
solutions are marginal, which indicates that the choice of high-order
stabilization for the limiting target
is of minor importance for this particular test problem.

\subsection{Two-dimensional solid body rotation}
To facilitate a 
comparison with the $\mathbb{P}_1/\mathbb{Q}_1$ version of
algebraic flux correction schemes and variational approaches
to shock capturing, let us
consider the solid body rotation benchmark \cite{convex,afc1,leveque}.
In this two-dimensional experiment, we solve the unsteady linear advection equation
\begin{equation}\label{sbr}
  \pd{u}{t}+\nabla\cdot(\mathbf{v}u)=0\quad\mbox{in}\quad\Omega=(0,1)^2
\end{equation}
using the divergence-free velocity field ${\bf v}(x,y)=2\pi(0.5-y,x-0.5)$ 
and the initial condition \cite{leveque}
\begin{subequations}\label{sbr_init}
\begin{equation}
  u_0(x,y)=\begin{cases}
  u_0^{\rm hump}(x,y)
  &\text{if}\  \sqrt{(x - 0.25)^2 + (y - 0.5)^2}\le 0.15, \\
  u_0^{\rm cone}(x,y)
  &\text{if}\ \sqrt{(x - 0.5)^2 + (y - 0.25)^2}\le 0.15, \\
  1 &\text{if}\ \begin{cases}
  \left(\sqrt{(x - 0.5)^2 + (y - 0.75)^2}\le 0.15 \right), \\
  \left(|x - 0.5| \ge 0.025,~ y\ge 0.85\right),
  \end{cases}\\
  0 &  \text{otherwise},
  \end{cases}
\end{equation}
where\vspace{-0.25cm}
\begin{align}
  u_0^{\rm hump}(x,y)&=    \frac14 + \frac14 \cos \left(
  \frac{\pi \sqrt{(x - 0.25)^2 + (y - 0.5)^2}}{0.15}\right),\\
  u_0^{\rm cone}(x,y)&= 1-\frac{\sqrt{(x - 0.5)^2 + (y - 0.25)^2}}{0.15}.
\end{align}
\end{subequations}
The so-defined initial data
undergoes counterclockwise rotation around the center
$(0.5,0.5)$ of the domain $\Omega$.
After each complete revolution (i.e., for $t\in\mathbb{N}$), the exact solution coincides with $u_0$.

\begin{figure}[!p]
  \centering\vspace{-0.75cm}
  \subfloat[HO-VMS-EV, $p=1$\label{fig:sbr_ho_vms_ev_p1}]{
    \begin{tabular}{c}
      {\scriptsize $||u_h-u_{\text{exact}}||_{L^1}=2.67\times 10^{-2}$} \\
      {\scriptsize $u_h\in[-4.6\times 10^{-3}, ~1.002]$} \\
    \includegraphics[scale=0.25]{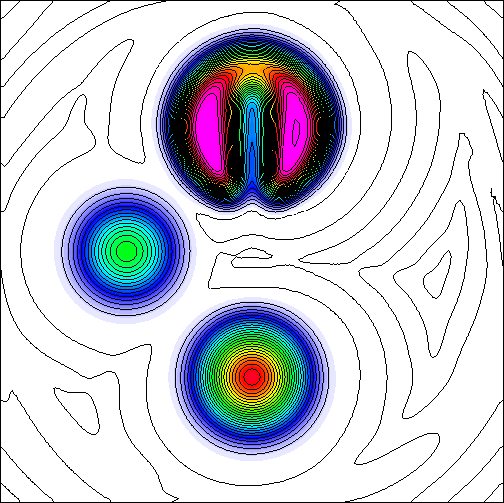}
    \end{tabular}
  }
  \subfloat[HO-VMS-EV, $p=2$\label{fig:sbr_ho_vms_ev_p2}]{
    \begin{tabular}{c}
      {\scriptsize $||u_h-u_{\text{exact}}||_{L^1}=1.63\times 10^{-2}$} \\
      {\scriptsize $u_h\in[-8\times 10^{-4}, ~0.9929]$} \\
      \includegraphics[scale=0.25]{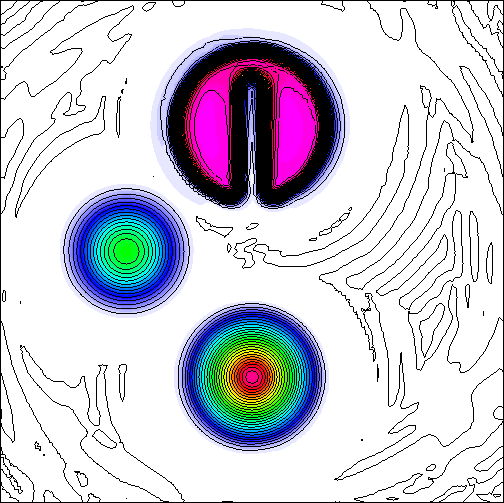}
    \end{tabular}
  }
  \subfloat[HO-VMS-EV, $p=4$\label{fig:sbr_ho_vms_ev_p4}]{
    \begin{tabular}{c}
      {\scriptsize $||u_h-u_{\text{exact}}||_{L^1}=9.96\times 10^{-3}$} \\
      {\scriptsize $u_h\in[-3.1\times 10^{-2}, ~1.044]$} \\
    \includegraphics[scale=0.25]{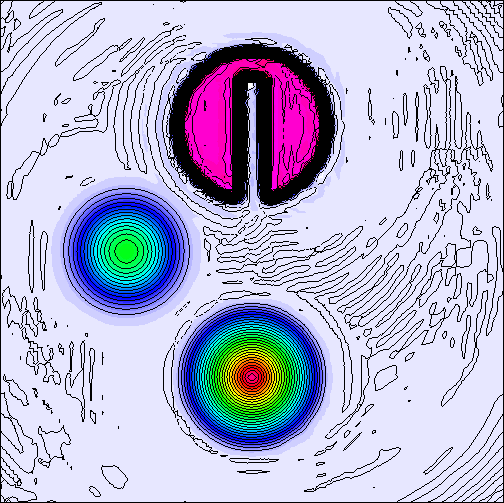}
  \end{tabular}
  }
  \vspace{5pt}

  \subfloat[HO-VMS-EV-BP, $p=1$\label{fig:sbr_ho_vms_ev_idp_p1}]{
    \begin{tabular}{c}
      {\scriptsize $||u_h-u_{\text{exact}}||_{L^1}=2.80\times 10^{-2}$} \\
      {\scriptsize $u_h\in[2.02\times 10^{-24}, ~0.9814]$} \\
    \includegraphics[scale=0.25]{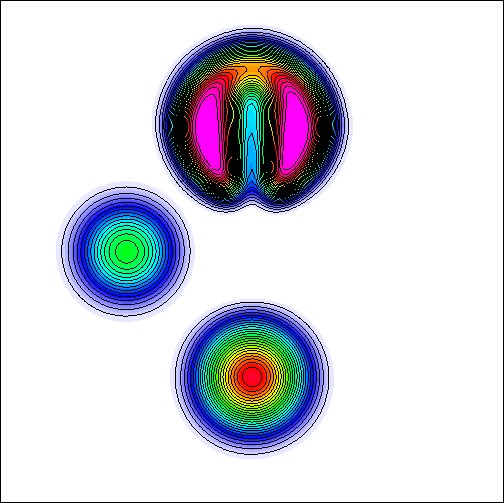}
    \end{tabular}
  }
  \subfloat[HO-VMS-EV-BP, $p=2$\label{fig:sbr_ho_vms_ev_idp_p2}]{
    \begin{tabular}{c}
      {\scriptsize $||u_h-u_{\text{exact}}||_{L^1}=2.49\times 10^{-2}$} \\
      {\scriptsize $u_h\in[6.06\times 10^{-23}, ~0.9865]$} \\
      \includegraphics[scale=0.25]{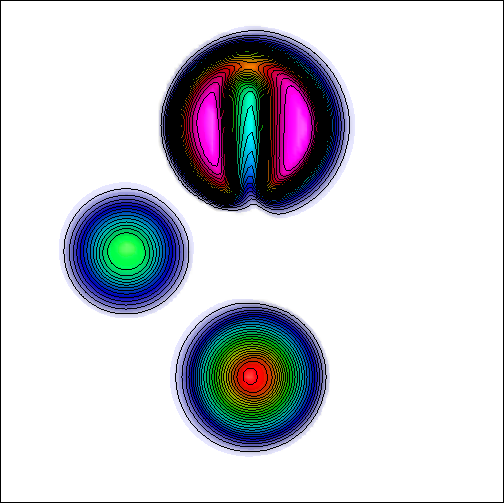}
    \end{tabular}
  }
  \subfloat[HO-VMS-EV-BP, $p=4$\label{fig:sbr_ho_vms_ev_idp_p4}]{
    \begin{tabular}{c}
      {\scriptsize $||u_h-u_{\text{exact}}||_{L^1}=2.11\times 10^{-2}$} \\
      {\scriptsize $u_h\in[2.59\times 10^{-19}, ~0.9862]$} \\
    \includegraphics[scale=0.25]{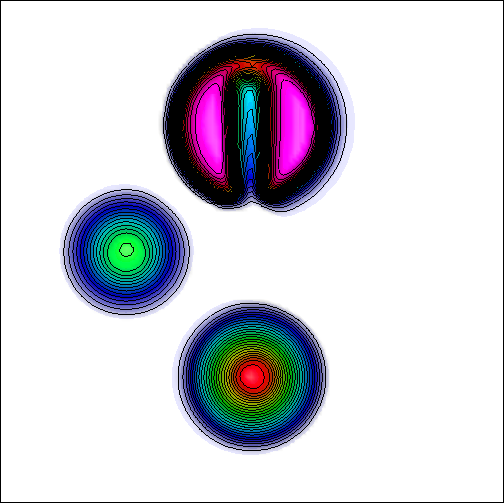}
  \end{tabular}
  }
  \vspace{5pt}

  \subfloat[HO-VMS-EV-FL, $p=1$\label{fig:sbr_ho_vms_ev_esLim_p1}]{
    \begin{tabular}{c}
      {\scriptsize $||u_h-u_{\text{exact}}||_{L^1}=3.67\times 10^{-2}$} \\
      {\scriptsize $u_h\in[2.16\times 10^{-11}, ~0.9195]$} \\
      \includegraphics[scale=0.25]{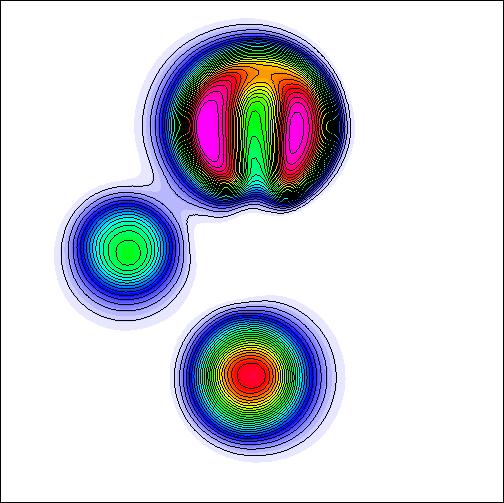}
    \end{tabular}
  }
  \subfloat[HO-VMS-EV-FL, $p=2$\label{fig:sbr_ho_vms_ev_esLim_p2}]{
    \begin{tabular}{c}
      {\scriptsize $||u_h-u_{\text{exact}}||_{L^1}=4.37\times 10^{-2}$} \\
      {\scriptsize $u_h\in[9.3\times 10^{-8}, ~0.8778]$} \\
      \includegraphics[scale=0.25]{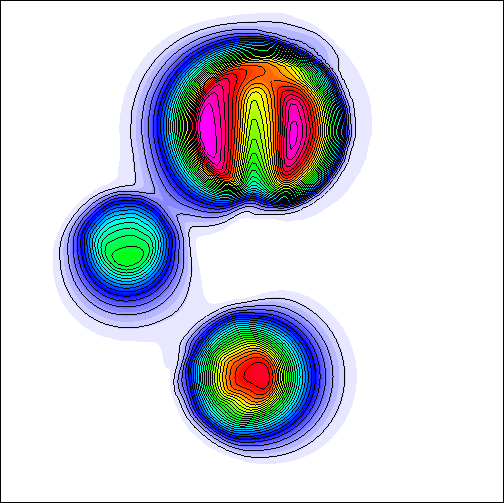}
    \end{tabular}
  }
  \subfloat[HO-VMS-EV-FL, $p=4$\label{fig:sbr_ho_vms_ev_esLim_p4}]{
    \begin{tabular}{c}
      {\scriptsize $||u_h-u_{\text{exact}}||_{L^1}=5.22\times 10^{-2}$} \\
      {\scriptsize $u_h\in[4.37\times 10^{-6}, ~0.7901]$} \\
      \includegraphics[scale=0.25]{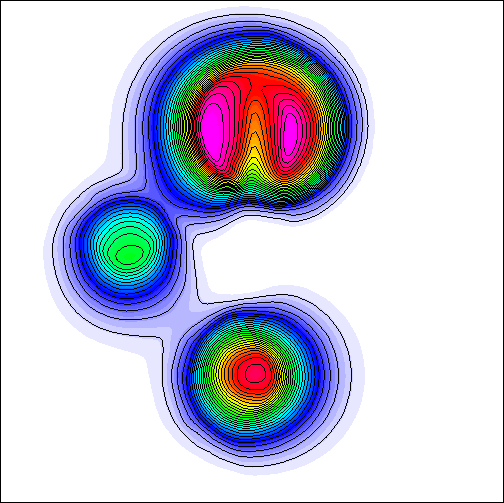}
    \end{tabular}
  }
  \caption{Solid body rotation problem \eqref{sbr} with initial condition \eqref{sbr_init}. Numerical solutions at $t=1$ obtained using $N_h=128^2$ DoFs.  In each diagram, we plot 30 contour lines corresponding to 
    a uniform subdivision of $\mathcal G = [0, 1]$.
  \label{fig:sbr}}
\end{figure}

Numerical solutions are evolved up to
$t=1$ using finite elements of degree $p=\{1,2,4\}$. For all values of $p$,
we choose the mesh size corresponding to $N_h=128^2$ DoFs. 
The results obtained with the entropy stable HO-VMS-EV method are shown in Fig. 
\ref{fig:sbr_ho_vms_ev_p1}-\ref{fig:sbr_ho_vms_ev_p4}. It can be seen
that the discontinuity-capturing effect of EV stabilization is not
enough to secure the IDP property w.r.t. $\mathcal G = [0, 1]$. The
flux-limited scheme HO-VMS-EV-BP yields the bound-preserving solutions
shown in Figs \ref{fig:sbr_ho_vms_ev_idp_p1}-\ref{fig:sbr_ho_vms_ev_idp_p4}.
The discontinuities are resolved in a nonoscillatory manner but the
imposition of local maximum principles results in unnecessary limiting at
smooth local extrema. To avoid the loss of high-order accuracy
around smooth traveling peaks, the local bounds of the subcell flux
limiting procedure can be relaxed using smoothness indicators, for
a presentation of which we refer the reader to \cite{Hajduk2020,convex2,CG-BFCT}. The results presented in Figs \ref{fig:sbr_ho_vms_ev_esLim_p1}-\ref{fig:sbr_ho_vms_ev_esLim_p4} were obtained with the entropy-aware HO-VMS-EV-FL version of the flux-limited scheme \eqref{fcorr}. In accordance with Remark \ref{remark:linear_adv}, the unnecessary imposition of the entropy stability condition \eqref{condES} increases the levels of numerical dissipation and the magnitude of the $L^1$~error. We conclude that limiting based on the local BP property
is sufficient for linear advection problems.

\subsection{Buckley-Leverett equation}
The first two-dimensional nonlinear problem that we consider is the 
Buckley-Leverett equation \cite{christov2008new,entropyCG}. The nonconvex
flux function of the nonlinear conservation law to be solved is  
\beq\label{bl}
\mathbf{f}(u)=
\frac{u^2}{u^2+(1-u)^2}
(1,1-5(1-u)^2).
\eeq
The computational domain is $\Omega_h=(-1.5,1.5)^2$. The
piecewise-constant initial condition is given by
\beq\label{bl_init}
u_0(x,y)=\begin{cases}
1 & \mbox{ if } x^2+y^2<0.5, \\
0 & \mbox{ otherwise}.
\end{cases}
\eeq
The exact solution of this nonlinear problem exhibits a rotating wave structure. For entropy stabilization purposes, 
we use $\eta(u)=\frac{u^2}2$. The corresponding entropy flux is
$\mathbf{q}(u)=(q_x(u),q_y(u))$, where
\begin{align}
q_x&=\frac14\left[
\frac{2(u-1)}{2u^2-2u+1}-\log(2u^2-2u+1)\right],\\
q_y&=\frac{1}{12}\left[
-20u^3+15u^2-\frac{9u+6}{2u^2-2u+1}-3\log(2u^2-2u+1)
-15\tan^{-1}(1-2u)\right].
\end{align}
An upper bound for the fastest wave speed can be found in \cite{christov2008new}. 
We overestimate it by using $\lambda=3.4$.

\begin{figure}[!h]
  \centering
  \subfloat[$p=1$]{
    \begin{tabular}{c}
      $u_h\in[0,~0.9993]$ \\ 
      \includegraphics[scale=0.28]{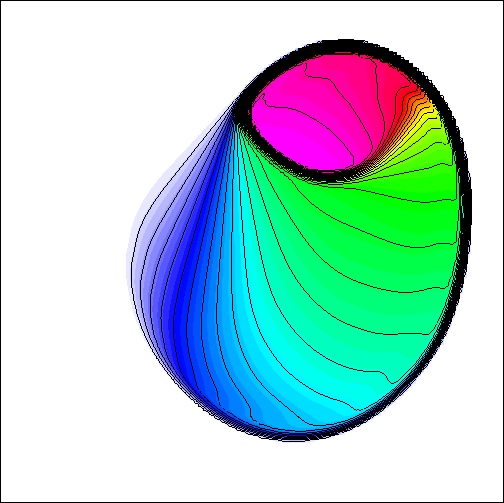}
    \end{tabular}
  }
  \subfloat[$p=2$]{
    \begin{tabular}{c}
      $u_h\in[0, ~0.993]$ \\ 
      \includegraphics[scale=0.28]{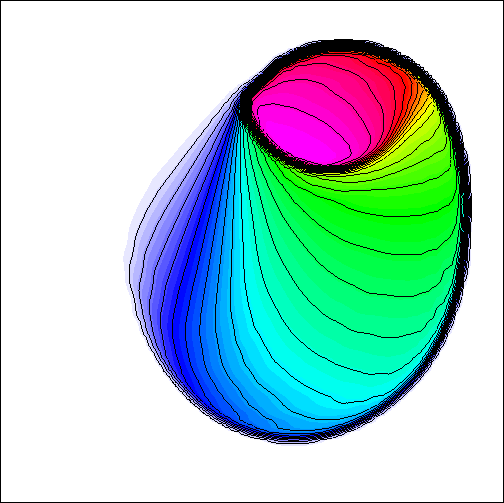}
    \end{tabular}
  }
  \subfloat[$p=4$]{
    \begin{tabular}{c}
      $u_h\in[0, ~0.9893]$ \\ 
      \includegraphics[scale=0.28]{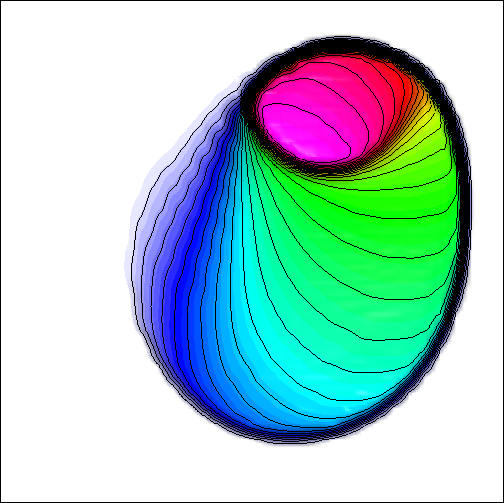}
    \end{tabular}
  }
  \caption{Buckley-Leverett equation \eqref{bl} with initial condition
    \eqref{bl_init}. Numerical solutions at $t=0.5$ obtained with HO-VMS-EV-FL
    using $N_h=128^2$ DoFs and Bernstein finite elements of degree $p=\{1,2,4\}$.
    In each diagram, we plot 30 contour lines corresponding to 
    a uniform subdivision of $\mathcal G = [0, 1]$.
  \label{fig:bl}}
\end{figure}

Simulations are performed using $N_h=128^2$ DoFs for
$p\in\{1,2,4\}$. The HO-VMS-EV-FL results at $t=0.5$ are shown in
Fig. \ref{fig:bl}.  They exhibit a crisp resolution of curved
shocks and are invariant domain preserving w.r.t. $\mathcal
G\in[0,1]$. The maximal values listed above the plots decrease
slightly as the polynomial degree $p$ is increased while keeping
$N_h$ fixed. However, the subcell flux limiting strategy makes
it possible to avoid a far more dramatic increase in the levels
of numerical dissipation due to extended stencils of high-order
finite element approximations (as reported, e.g., in \cite{CG-BFCT}).

\subsection{KPP problem}
In the last numerical example, we consider the KPP problem \cite{Guermond2016,Guermond2017,kpp}, a challenging nonlinear test for verification of entropy stability properties. Equation \eqref{ibvp-pde} with 
the nonconvex flux function
\beq
\mathbf{f}(u)=(\sin(u),\cos(u))
\eeq
is solved in the computational domain
$\Omega_h=(-2,2)\times(-2.5,1.5)$ using the initial
condition
\beq \label{kpp_init}
u_0(x,y)=\begin{cases}
\frac{7\pi}{2} & \mbox{if}\quad \sqrt{x^2+y^2}\le 1,\\
\frac{\pi}{4} & \mbox{otherwise}.
\end{cases}
\eeq
The entropy flux corresponding to $\eta(u)=\frac{u^2}2$ is
$\mathbf{q}(u)=(u\sin(u)+\cos(u),u\cos(u)-\sin(u))$.
A simple upper bound for the 
maximal speed is $\lambda=1$. More accurate estimates
can be found in \cite{Guermond2017}.

\begin{figure}[!h]
  \centering
  \subfloat[$p=1$]{
    \begin{tabular}{c}
      $u_h\in[0.7854, ~10.9955]$ \\ 
      \includegraphics[scale=0.28]{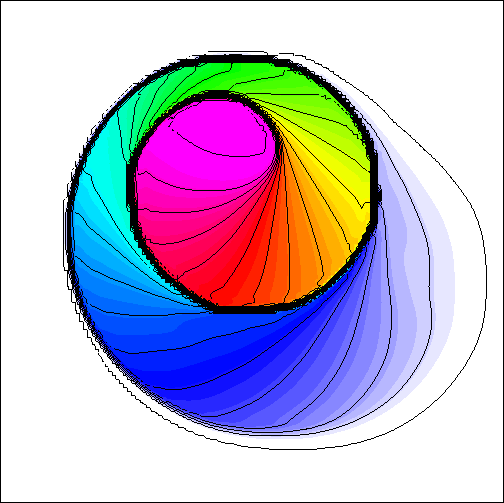}
    \end{tabular}
  }
  \subfloat[$p=2$]{
    \begin{tabular}{c}
      $u_h\in[0.7854, ~10.9955]$ \\ 
    \includegraphics[scale=0.28]{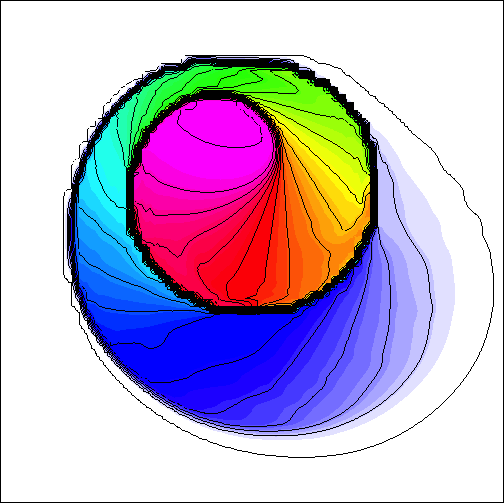}
  \end{tabular}
  }
  \subfloat[$p=4$]{
    \begin{tabular}{c}
      $u_h\in[0.7854, ~10.9955]$ \\ 
      \includegraphics[scale=0.28]{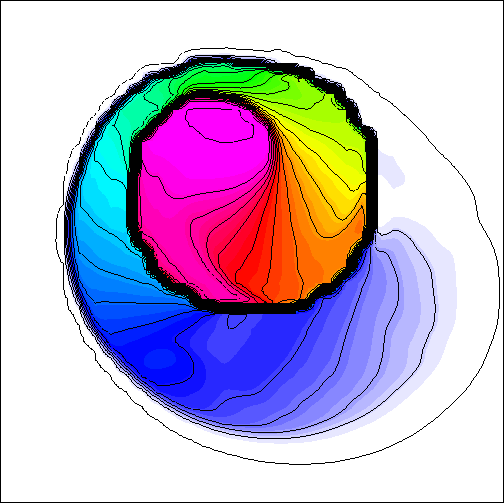}
    \end{tabular}
  }
  \caption{KPP problem \cite{kpp} with initial condition \eqref{kpp_init}.
    Numerical solutions at $t=1$ obtained with HO-VMS-EV-FL
    using $N_h=128^2$ DoFs and Bernstein finite elements of degree $p=\{1,2,4\}$.
    In each diagram, we plot 30 contour lines corresponding to 
    a uniform subdivision of $\mathcal G = [\frac14\pi, ~\frac72\pi]$.
  \label{fig:kpp}}
\end{figure}

Similarly to the Buckley-Leverett problem, the entropy solution of 
the KPP problem exhibits a two-dimensional rotating wave structure.
The main challenge of this test is to
prevent possible convergence to wrong weak solutions. Even
bound-preserving high-resolution
schemes may fail to preserve the thin gap
between the twisted shocks if no entropy viscosity is added
\cite{Guermond2017,entropyCG}. The results displayed in Fig.~\ref{fig:kpp}
were obtained with HO-VMS-EV-FL using 
$N_h=128^2$ DoFs for Bernstein finite elements
of degree $p=\{1,2,4\}$. The snapshots correspond to the final time $t=1$
and reproduce the rotating wave structure of the entropy solution
correctly (cf. \cite{Guermond2017,entropyDG,entropyCG}). The two shocks remain
clearly separated and the ranges of the numerical solutions stay
in the invariant set $\mathcal G=[\frac14\pi, ~\frac72\pi]$. 

\section{Conclusions}

The presented research was aimed at exploring the aspects of entropy stability in the context of high-order continuous finite element approximations to hyperbolic conservation laws. We proved that the continuous Galerkin method is square entropy conservative, endowed it with high-order stabilization terms, and designed property-preserving limiters for high-order Bernstein finite elements. It is hoped that the proposed methodology paves the way for further analysis and design of nonlinear high-resolution finite element schemes equipped with entropy correction procedures. In particular, we envisage that extensions of the new entropy fixes to hyperbolic systems and discontinuous Galerkin methods should be relatively straightforward. More challenging open problems include theoretical investigations of the steady-state limit, development of efficient iterative solvers for nonlinear discrete problems,  and provable preservation of entropy stability in fully discrete flux-corrected schemes.

\medskip
\paragraph{\bf Acknowledgments}

The work of Dmitri Kuzmin was supported by the German Research Association (DFG) under grant KU 1530/23-1. The authors would like to thank Hennes Hajduk (TU Dortmund University) for suggesting an improved version of the subcell flux decomposition.


\end{document}